\documentclass[12pt]{article}
\usepackage{amssymb}
\usepackage[centertags]{amsmath}
\usepackage{pstricks}
\usepackage{pst-node}
\usepackage{latexsym}

\newtheorem{Lemma}{Lemma}
\newtheorem{Proposition}[Lemma]{Proposition}

\newtheorem{Theorem}[Lemma]{Theorem}

\newtheorem{Corollary}[Lemma]{Corollary}

\newcommand{\bck}{\!\!\!}
\newcommand{\ed}{\ \stackrel{d}{=} \ }

\newcommand{\GG}{\mbox{${\cal G}$}}
\newcommand{\FF}{\mbox{${\cal F}$}}
\newcommand{\PP}{\mbox{${\cal P}$}}

\newcommand{\eps}{\varepsilon}

\newcommand{\bN}{{\bf N}}
\newcommand{\bR}{{\bf R}}
\newcommand{\bZ}{{\bf Z}}

\newcommand{\la}{{\Lambda}}
\newcommand{\mg}{{\upsilon}}
\def\endpf{$\Box$}

\begin{document}

\date{\today}
\title{The spatial $\la$-coalescent}
\author{Vlada Limic\thanks{Research supported in part by an NSERC
research grant.}\\
Department of Mathematics\\
The University of British Columbia\\
Vancouver, BC V6T 1Z2
\and
Anja Sturm\\ 
Department of Mathematical Sciences\\
University of Delaware\\
Newark, DE 19711}
\maketitle

\abstract{
This paper extends the notion of the $\la$-coalescent of Pitman (1999)
to the spatial setting.  The partition elements of the spatial
$\la$-coalescent
migrate in a (finite) geographical space  and may only coalesce if
located at the same site of the space. We characterize the $\la$-coalescents
that come down from infinity, in an analogous way to Schweinsberg (2000).
Surprisingly, all spatial coalescents that come down from infinity, also come
down from infinity in a uniform way. This enables us to study space-time
asymptotics of spatial $\la$-coalescents on large tori in $d\ge 3$
dimensions.
Our results generalize and strengthen those of Greven et al.~(2005), who
studied
the spatial Kingman coalescent in this context.
}
        \vspace{\fill}

\par
        \emph{AMS 2000 Subject Classification.} Primary 60J25,
        60K35 \\
\par
        \emph{Key words and phrases.}
coalescent, $\la$-coalescent, structured coalescent, limit theorems
          \thispagestyle{empty}

\par
         
\newpage

\section{Introduction}

The {\em $\la$-coalescent}, sometimes also called the coalescent with multiple collisions,
 is a Markov process $\Pi$ whose state space is the set of partitions of the positive integers. 
The
{\em standard $\la$-coalescent} $\Pi$ starts at the partition of the
positive integers into singletons, and its restriction to
$[n]:=\{1,\dots,n\},$ denoted by $\Pi_{n},$ is the $\la$-coalescent
starting with $n$ initial partition elements. The measure $\la,$
which is a finite measure on $[0,1],$ dictates the rate of coalescence
events, as well as how many of the (exchangeable) partition elements, which we will also 
refer to as {\em blocks},
may coalesce into one at any such event.
The $\la$-coalescent was introduced by Pitman \cite{jP99},
and also studied by Schweinsberg \cite{jS00a}.
It was obtained as a limit
of genealogical trees in a Moran-like model by
Sagitov \cite{sS99}.

The well-known Kingman coalescent
\cite{Kin82} corresponds to the $\la$-coalescent with
$\la(dx)=\delta_0(dx)$, the unit atomic measure at $0$.
For this coalescent, each pair of current partition elements coalesces
at unit rate, independently from other pairs.
Papers \cite{aldsurvey} and \cite{EP98} are devoted to
stochastic coalescents where again only pairs of partitions are allowed to
coalesce, but the coalescence rate is not uniform over all pairs.
The survey \cite{aldsurvey} gives many pointers to the literature.
The $\la$-coalescent generalizes the Kingman coalescent in the sense
that now any number of  partition elements may merge into one at a
coalescence
event, but the rate of coalescence for any $k$-tuple of partition elements
depends still only on $k.$ The first example of such a $\la$-coalescent
(other than the Kingman coalescent) was studied by  Bolthausen and
Sznitman \cite{BS98}, who were interested in the special case where
$\la(dx)$ is Lebesgue measure on $[0,1]$ in connection with spin glasses.
Bertoin and Le Gall \cite{BLG00} observed a correspondence of this particular
coalescent to the genealogy of continuous state branching processes (CSBP).
More recently, Birkner et al.~\cite{BBCEMSW05} extended
this correspondence to stable CSBP's to  $\la$-coalescents, where $\la$ is given
by a Beta-distribution.
Berestycki et al.~\cite{BBS} 
use this correspondence to study fine
small time properties of the corresponding coalescents.

A further generalization of the $\la$-coalescents,
known as the {\em coalescents with simultaneous multiple collisions},
was originally studied by
M\"ohle and Sagitov \cite{MS01} and
Schweinsberg \cite{jS03}.
Further connections to bridge processes and generalized
Fleming-Viot processes were discovered by
Bertoin and Le Gall \cite{BLG03},
and to asymptotics of genealogies during
selective sweeps, by
Durrett and Schweinsberg \cite{DS05}.

Our first goal, in Section \ref{CON}, is to extend the notion of the $\la$-coalescent to the spatial
setting. Here, partition elements migrate in a geographical
space and may only coalesce while sharing the same location.
Earlier works on variants of spatial coalescents,
sometimes
also referred to as structured coalescents, have all assumed
Kingman coalescent-like behavior, and include
Notohara (1990) \cite{mN90}, Herbots (1997)
\cite{hH97},
and more recently Barton et al.~\cite{BES04} in the case of finite initial
configurations, and Greven et al.~\cite{GLW} with infinite initial
states.
A related model has been studied by Z\"ahle et al.~\cite{ZCD}
on two-dimensional tori.

\smallskip
In most of this paper we assume that $\la$ is
a finite measure on $[0,1]$ without an atom at $0$ or
at $1$, such that $\la([0,1])>0$.
At the end of Section \ref{CON} we
comment on how atoms at $0$ or $1$ would change the
behavior of the coalescent.

Define for $2\leq k\leq b$, $k,b$ integers,
\begin{equation}
\label{Edeflambk}
\lambda_{b,k}:=\int_{[0,1]} x^{k-2}(1-x)^{b-k} d\la(x).
\end{equation}
The parameter $\lambda_{b,k}\geq 0$ is the rate at which $k$ blocks
coalesce when the current configuration has $b$ blocks.
Extend the
definition by setting $\lambda_{b,k}=0$ for $b=1$ or $b=0$, $k\in \bN$.
Define in addition
\begin{equation}
\label{Edeflam}
\lambda_b:=
\sum_{k=2}^b{b \choose k}\lambda_{b,k},
\end{equation}
and
\begin{equation}
\label{Edefgam}
\gamma_b:=
\sum_{k=2}^b {b \choose k}(k-1)\lambda_{b,k}.
\end{equation}
Note that $\lambda_b$ is the total rate of coalescence when
the configuration has $b$ blocks, and
that $\gamma_b$
is the total rate of decrease in the number of blocks when
the configuration has $b$ blocks. From the above definitions, one may 
already observe (see also proof of Theorem \ref{Tconstr})
that the $\la$-coalescent
can be derived from a Poisson point process on $\bR_+ \times [0,1]$ ($\bR_+:=[0,\infty)$) with intensity
measure $ dt x^{-2} d\la(x):$ If $(t,x)$ is an atom of this Poisson point process, then at 
time $t,$ we mark each block independently with probability $x,$ and subsequently merge
all marked blocks into one. 

Now consider a finite graph $\GG$, and denote by $|\GG|$ the number of its
vertices.
Call the vertices of $\GG$ {\em sites}. 
Consider a process started from a finite configuration
of $n$ blocks on sites
in ${\cal G}$ where we allow only two types of transitions,
referred to as {\em coalescence} and {\em migration} respectively:
\begin{itemize}
\item[(i)] at each site blocks coalesce according to the $\la$-coalescent,
\item[(ii)] the location process
of each block is an independent continuous Markov chain on $\GG$
with jump rate $1$ and transition probabilities $p(g_i,g_j), g_i,g_j \in \GG$.
\end{itemize}
The original $\la$-coalescent of \cite{jP99} and \cite{jS00a}
corresponds to the setting where $|\GG|=1$, so migrations are impossible.
The spatial $\la$-coalescent started from a finite configuration
$\{(1,i_1),\ldots,(n,i_n)\}$
is a well-defined strong Markov process (chain)
with state space being the set
of all partitions of $[n]=\{1,\dots,n\}$ labeled by their location in
${\cal G}$.
This will be stated precisely in
Theorem \ref{Tconstr} of Section \ref{CON} which is devoted to
the construction of spatial $\la$-coalescents $\Pi^{\ell}$ with general (possibly infinite)
initial states.

After constructing the general spatial $\la$-coalescent, we turn to characterizing those that come down from infinity in Section \ref{CDI}.
Schweinsberg \cite{jS00a} shows that if
\begin{equation}
\label{Edown}
\sum_{b \geq 2} \frac{1}{\gamma_b}<\infty
\end{equation}
holds, then the (non-spatial)
$\la$-coalescent started with infinitely many blocks at time $0$
immediately
comes down from infinity, that is, the number of its blocks
at all times $t>0$ is finite with
probability 1; otherwise, the $\la$-coalescent stays
infinite forever, meaning that it contains infinitely many blocks at all
times
$t>0$ with probability 1.

The goal of Section \ref{CDI}
is to show that the spatial $\la$-coalescent
inherits this property of either coming down from
infinity or staying infinite, from its nonspatial counterpart.
More precisely, let $(\Pi^\ell(t))_{t\geq 0}$ be the $\la$-coalescent
constructed in
Theorem \ref{Tconstr}, and
denote by $\#\Pi(t)$ its size
at time $t$, i.e. the total number of blocks in $\Pi^\ell(t)$, with any
label.
In Lemma \ref{LETn} and Proposition \ref{PzoE}
we show that condition (\ref{Edown}) implies $P[\#\Pi(t)<\infty, \forall t>0]=1$, even if
the initial configuration $\Pi(0)$ contains infinitely many blocks.
In this case we say
that the spatial $\la$-coalescent {\em comes down from infinity}.
In Proposition \ref{PzoE} we also show via a coupling to the
non-spatial coalescent
that if (\ref{Edown}) does not hold, provided $\#\Pi(0)=\infty$ and
$\la$ has no atom at $1$, then $P[\#\Pi(t)=\infty, \forall t>0]=1$.
In this case we say
that the spatial $\la$-coalescent {\em stays infinite}.
We note here that 
the statement of Lemma \ref{LETn} (saying that $\sup_{n}E[T_n]<\infty,$ where
$T_n$ is the time until there are on average two 
blocks per site if there are initially $n$ blocks per site)  
extends to the spatial coalescent for which the migration mechanism
may be more general, for example non-exponential or depending 
on the coalescence mechanism.

In Section \ref{UNI} we continue the study of the time
$T_n.$ In particular, in Theorem   \ref{Tuniform} we obtain an upper bound
on its expectation that is
not only uniform in $n$ but also, somewhat surprisingly,
in the structure (size) of $\GG$. In this case, we say that
the coalescent {\em comes down from infinity uniformly.}
The argument of Theorem \ref{Tuniform} relies on the
independence of the coalescence and migration mechanisms.

Our final goal, in Section \ref{UNI}, is to study space-time asymptotic properties
of $\la$-coalescents that  come down from infinity uniformly 
on large finite tori at time scales on the order of the volume.
In \cite{GLW}, this asymptotic behavior was studied for
the spatial Kingman coalescent where
$\la=\gamma \delta_0$ for some $\gamma>0$.
It is interesting that on appropriate space-time scales,
the scaling limit is again (as in \cite{GLW}) the
Kingman coalescent, with only its starting configuration depending
on the specific properties of the underlying $\la$-coalescent.
We obtain
functional limit theorems for the partition structure and
for the number of partitions, 
in Theorems \ref{TnpartCon} and \ref{TuniCon} respectively.

\section{Construction of the coalescent}
\label{CON}
The construction of the spatial coalescent
on an appropriate state space follows quite standard steps. 
The construction below is inspired by those in  Evans and Pitman \cite{EP98},
Pitman \cite{jP99}, and Berestycki \cite{jB04}. 

Let $\PP$ be the set of partitions on $\bN,$ which can be identified
with the set of equivalence relations on $\bN.$ 
Any $\pi \in \PP$ can be represented uniquely by
$\pi=(A_{1},A_{2},A_{3},\dots)$ where $A_{j}\subset \bN$ for
$j \geq 1$ are called the the {\em blocks}  of $\pi$, indexed according to
the increasing ordering of the set $\{\min A_j: j \geq 1\}$ that contains the
smallest element of each block.
So in particular $\min A_{n-1} < \min A_n$, for any $n\geq 2$. 
Likewise, we define for any $n \in \bN,$ 
$\PP_n$ as the set of partitions of $[n],$  and for 
$\pi \in \PP_n$ we have $\pi=(A_{1},A_{2},\dots, A_{n})$ in an 
analogous way.
We will write $A \in \pi$ if $A\subset \bN$ is a block of $\pi,$ and 
\[
A_i \sim_{\pi} A_j
\] 
if $A_i,A_j \subset \bN$ and $A_{i} \cup A_{j} \subset A$   
for some (unique) $A \in \pi.$
If the number of blocks of $\pi,$ denoted by $\#\pi,$ is
finite, then set $A_{j}=\emptyset$ for all $i>\#\pi$.

For concreteness in the rest of the paper, let $|\GG|=\mg$ for $\mg$ a positive
integer  and let the vertices of $\GG$ be $\{g_1,\ldots,g_\mg\}.$
The spatial coalescent takes values in the set $\PP^{\ell}$ of 
partitions on $\bN$, indexed as described above, and labelled by $\GG$, so
\[
\PP^{\ell}:=\{(A_j,\zeta_j) : A_j \in \pi, \zeta_j \in \GG, \pi \in \PP, j\geq 1\}.
\]
Similarly, the coalescent started from $n$ blocks takes values
in $\PP^{\ell}_n:= 
\{(A_j,\zeta_j): A_j  \in \pi, 
\zeta_j \in \GG, \pi \in \PP_n, 1\leq j\leq n\}.$
Here, the $\zeta_j \in \GG$ is the {\em label} (or location)
of $A_j\in \pi$, $j\geq 1$. 
Set $\zeta_j=\partial \not\in \GG$ if
$A_{j}=\emptyset$. 
For any element $\pi \in
\PP^{\ell}$ or $\pi \in \PP^{\ell}_{n}$ with $n\geq m$ 
define $\pi|_m \in \PP^{\ell}_m$ as the labeled partition
induced by $\pi$ on $\PP^{\ell}_m$.
We equip $\PP^{\ell}$ with the
metric
\begin{equation}
\label{Emetric}
d(\pi,\pi')= \sup_{m\in \bN} 2^{-m}1_{\{\pi|_m \neq \pi'|_m\}}; \quad
\pi, \pi' \in \PP^{\ell},
\end{equation}
and likewise $\PP^{\ell}_n$ with the metric
\begin{equation*}
d_n(\pi,\pi')= \sup_{m \leq n} 2^{-m}1_{\{\pi|_m \neq \pi'|_m\}}; \quad
\pi, \pi' \in \PP^{\ell}_n.
\end{equation*}
It is easy to see that
$(\PP^{\ell}_n,d_n)$ and $(\PP^{\ell},d)$ are \emph{both} 
compact metric spaces, and that $d(\pi,\pi') = \sup_n d_n(\pi|_n,\pi'|_n)$.

Note that $\PP^{\ell}$ can be interpreted as a subspace
of the infinite product space 
$\PP^{\ell}_{\prod}:=
(\PP^{\ell}_1,\PP^{\ell}_2,\PP^{\ell}_3,\dots)$ 
endowed with the metric 
$d(\pi,\pi')=\sup_{n} d_{n}(\pi_n,\pi'_n)$ for 
$\pi=(\pi_1,\pi_2,\dots),$ $\pi'=(\pi'_1,\pi'_2,\dots) \in \PP^{\ell}_{\prod}$,
by identifying
$\pi \in \PP^{\ell}$ with $(\pi|_1,\pi|_2,\pi|_3,\dots)$.
Note that
this metric induces the product topology on $\PP^{\ell}_{\prod}$ and
that an element $(\pi_1,\pi_2,\dots)$ of $\PP^{\ell}_{\prod}$ is
also an element of $\PP^{\ell}$ if and only if it fulfils the
following \emph{consistency relationship,}
\begin{equation}
\label{consistency}
\pi_{n+1}|_n =\pi_{n} \quad \mbox{for all }n \in \bN.
\end{equation}

In the rest of the paper, whenever $\Pi^{\ell}$
is a spatial coalescent process, we denote by $\Pi$ the partition (without
the labels of the blocks) of $\Pi^\ell$, and by 
\[
(\#\Pi(t))_{t\geq 0}
\]
the corresponding total number of blocks process.
Thus $\#\Pi(t)$ is the number (finite or infinite) of blocks in
$\Pi(t)$, or equivalently, in $\Pi^\ell(t)$.

With the above notation we are finally able to construct the 
spatial $\la$-coalescent started from potentially infinitely many blocks, as stated
in the following theorem.
Recall the migration mechanism stated in the introduction: 
each block performs an independent continuous Markov chain
on $\GG$ with jump rate $1$ and transition kernel $p(\cdot,\cdot)$.
\begin{Theorem}
\label{Tconstr}
Assume that $\la$ has no atom at $0$.
Let $\GG$ be a finite graph with vertex set 
$\{g_1,\ldots,g_\mg\}$. 
Then, for each $\pi \in \PP^{\ell},$ there exists a 
c\`adl\`ag Feller and strong Markov process $\Pi^{\ell}$
on $\PP^{\ell}$, called the {\emph spatial $\la$-coalescent}, such that 
$\Pi^{\ell}(0)=\pi$ and 
\begin{itemize}
\item[(i)] blocks with the same label coalesce according to a (non-spatial)
$\la$-coalescent,
\item[(ii)] each block of label $g_i\in \GG$ changes its label to $g_j\in \GG$
at rate $p(g_i,g_j)$ as mentioned in introduction.
\end{itemize}
This process also satisfies
\begin{itemize}
\item[(iii)]
$(\Pi^{\ell}(t)|_n)_{ t\geq 0}$ is a 
spatial $\la$-coalescent started from $\Pi^{\ell}(0)|_n$,
\end{itemize}
and its law is characterized by (iii) and the initial configuration $\pi$.
\end{Theorem}

{\em Proof.}
In order  to define a
c\`adl\`ag Markov process $\Pi^{\ell}$ with values in $\PP^{\ell}_{\prod}$ such 
that $\Pi^{\ell}_n:=\Pi^{\ell}|_n$ is a spatial coalescent
starting 
at $\Pi^{\ell}(0)|_n \in
\PP^{\ell}_n$ for any $\Pi^{\ell}(0) \in \PP^{\ell},$
we will make use of suitably chosen Poisson point processes. 

For each $i \in [\mg]$ let $N_{i}$ be an independent Poisson point process on 
$\bR_{+} \times [0,1] \times \{0,1\}^{\bN}$ with intensity measure 
$dt \: x^{-2} \la(dx) P_{x}(d\xi),$ where 
$\xi=(\xi_1,\xi_2,\dots)$ is a random vector whose entries $\xi_j$ are i.i.d.~Bernoulli($x$)  under $P_x$, defined on some probability space 
$(\Omega,\FF,\FF_t, P).$  

Let $\delta_n$ denote the Kronecker delta measure with unit atom at $n$.
Let $M$ be another independent Poisson point process on the same probability
space $\Omega$ with values in $\bR_{+} \times \bN \times \GG^{\mg}$ 
and intensity measure given by 
$dt \sum_{k=1}^{\infty} \delta_{k}(dm) 
P^{\mg}(ds_1,\ldots,ds_{\mg}),$
where $P^{\mg}$ is the joint law of $\mg$ independent $\GG$-valued
random variables $S_1,\ldots,S_\mg$, such that $P(S_{g_i}=g_j)=p(g_i,g_j)$,
$g_i,g_j\in \GG$.

Using the above random objects define a spatial $\la$-coalescent
with $n$ initial blocks, 
$\Pi^{\ell}_n,$ on $\Omega$ for each $n \geq 1$ as follows: 
At any atom $(t,x,\xi)$ of $N_{i}$,
all blocks $A_j(t-)$ with $\zeta_j(t-)=g_i$ and $\xi_j=1$
coalesce together into a new labeled block
$(\bigcup_{j,\xi_j=1,\zeta_j(t-)=i}  A_j(t-), g_i)$;
at any atom $(t,m,(s_1,\ldots,s_{\mg}))$ of $M$ 
we set $\zeta_m(t)=s_{\zeta_m(t-)}$ provided $m \leq \#\Pi_n(t-)$,
otherwise nothing changes.
For all other $t \geq 0$ we set $\Pi^{\ell}_n(t)=\Pi^{\ell}_n(t-).$
Note that coalescence causes immediate reindexing 
(or reordering) of blocks
that have neither participated in  coalescence nor in migration, and that 
this reindexing operation decreases each index by a non-negative amount.

Since the sum of the above defined jump rates of $\Pi^{\ell}_n$ is finite
it follows immediately that $\Pi^{\ell}_n$ is a well defined c\`adl\`ag
Markov process on $\Omega$ for each $n \geq 1$
 therefore inducing a 
c\`adl\`ag Markov process $\Pi^{\ell}$ on $\PP^{\ell}_{\prod}$.
It is important to note that each $\Pi^{\ell}_n$ so constructed is 
a $\la$-coalescent started from $\Pi^{\ell}(0)|_{n}.$ 
Since $\Pi^{\ell}_{n+1}(0)|_n = \Pi^{\ell}_{n}(0)$ and since 
clearly the consistency condition
(\ref{consistency}) is preserved under each transition 
of $\Pi^{\ell}_{n+1}$ in the construction
(this is not always a transition for $\Pi^{\ell}_{n}$), we
have $\Pi^{\ell}_{n+1}(t)|_n = \Pi^{\ell}_{n}(t)$ for all $t \geq 0.$ 
Therefore, $(\Pi^{\ell}(t))_{t \geq 0}$ constructed by $\Pi^{\ell}(t)|_n:=\Pi^{\ell}_n(t), n\geq 1,t\geq 0$ 
is well-defined.
It follows that $\Pi^{\ell}$ is a c\`adl\`ag Markov process with 
values in $\PP^{\ell},$ which clearly satisfies 
properties (i)-(iii), and uniqueness in distribution 
follows similarly.
 
In order to verify that 
the semigroup 
$T_t \varphi(\pi) := E[\varphi(\Pi^{\ell}(t))| \Pi^{\ell}(0)=\pi]$  
is a Feller-Dynkin semigroup it now suffices to check the
following two properties (see \cite{RW1} III (6.5)-(6.7)):
(i) For any continuous (bounded) real valued 
function $\varphi$ on $\PP^{\ell}$ and all 
$\pi \in \PP^{\ell}$ we have
\begin{equation*}  
\lim_{t \rightarrow 0+} T_t \varphi(\pi) =\varphi(\pi),
\end{equation*} 
and 
(ii)
for any continuous (bounded)
real valued function $\varphi$ on $\PP^{\ell}$
and all $t>0,$ $\pi \mapsto T_t\varphi(\pi)$ is continuous (and bounded).

Note that (i) is an immediate consequence of the right-continuity of the
paths and continuity with respect to (\ref{Emetric}).
One can easily argue for (ii):
if $\Pi^{\ell,k}$ is the spatial coalescent started from $\pi^k$
and $\Pi^{\ell}$ is the spatial coalescent started from $\pi$ 
such that 
$\lim_{k \rightarrow \infty}
\pi^k = \pi \in \PP^{\ell}$,
then, due to the definition of the metric (\ref{Emetric}) on $\PP^\ell$,
there exists for all $k \in \bN$ an $m=m(k)$ such that 
$\pi^{k}|_{m} =\pi|_m$, with the property
$m(k)\to \infty$ as $k\to \infty$.
This implies that
one can construct a coupling of $\Pi^{\ell,k}$ and $\Pi^{\ell}$ (using the same
Poisson point processes for all) such that
$\Pi^{\ell,k}(t)|_m =\Pi^{\ell}(t)|_m$ for all $t\geq 0$.
Hence $d(\Pi^{\ell,k}(t),\Pi^{\ell}(t))\leq 2^{-(m(k)+1)}$ for all $t\geq 0$ 
and, since $m(k)\to \infty$, we conclude that the second property
holds due to the continuity of $\varphi$.
Given that $T_t$ is a Feller-Dynkin semigroup the strong Markov property holds.
\hfill \endpf

\smallskip
{\em Remark.}
A variation of the above construction could be repeated for 
the cases where $\la$ has an atom at $0$.
This would correspond to superimposing Kingman coalescent type transitions
on top of the Poisson process induced coalescent events.
One easily observes that all such coalescents come down from infinity.
Also note that an atom of $\Lambda$ at $1$ implies complete collapse 
in finite time, even if the
coalescent corresponding to the measure $\Lambda(\cdot \cap [0,1))$ stays infinite.
See \cite{jP99} for further discussion of atoms.
\hfill $\diamond$

\smallskip
{\em Remark.}
We stated Theorem \ref{Tconstr} for $|\GG|<\infty.$ The case $|\GG|=\infty$ needs a little more work if we also want to be able to start with an infinite configuration $\pi \in \PP^{\ell}.$ However, for $\pi \in \PP_n^{\ell},$ a finite starting configuration, the Poisson point process construction in the proof of the theorem 
immediately yields the desired process.
This fact will be useful in Section \ref{ASY} where we consider 
$\GG=\bZ^d.$ 
\hfill $\diamond$

\section{Coming down from infinity}
\label{CDI}
In this section, we first obtain
estimates on the coalescence rates 
and the rates of decrease in the number of blocks, both  in the non-spatial and the spatial 
setting. 
Several of these estimates will be applied to
showing that the spatial $\la$-coalescent comes down 
from infinity if and only if (\ref{Edown}) holds.

It is easy to see, using definitions (\ref{Edeflambk})-(\ref{Edefgam}), that
\begin{equation}
\label{Elam_gam}
\lambda_b= \int_{[0,1]} \frac{1-(1-x)^b-bx(1-x)^{b-1}}{x^2}d\la(x) \ , \
\gamma_b= \int_{[0,1]} \frac{b x-1+(1-x)^b}{x^2}d\la(x).
\end{equation}
The following lemma is listing some 
facts, which are based on (\ref{Elam_gam}) and some simple computations.
\begin{Lemma}
\label{Linit}
We have the following estimates:
\begin{itemize}
\item[(i)]$\lambda_{b+1}-\lambda_b=\int_{[0,1]} b(1-x)^{b-1} d\la(x)$ 
for $b\geq 2$, in particular $\lambda_b\leq \lambda_{b+1}\leq
3\lambda_b$,
\item[(ii)] $\gamma_{b+1}-\gamma_b=\int_{[0,1]} (1-(1-x)^{b}) x^{-1} d\la(x)\geq
0$.
\end{itemize}
\end{Lemma}
{\em Proof.}
(i) Note that
$-(1-x)^{b+1}-(b+1) x (1-x)^{b}+(1-x)^{b}+b x (1-x)^{b-1}=$
$(1-x)^{b-1}(-(1-x)^2 -(b+1)x(1-x)+(1-x)+bx)$,
and that the term in the parentheses equals $b x^2$.
Combined with (\ref{Elam_gam}), this  gives the initial statement of the lemma.
The first inequality $\lambda_{b+1}\geq \lambda_b$ is immediate.
The second inequality follows again from (\ref{Elam_gam}), 
by integrating the following inequality with respect to $\la$ 
\[
b  (1-x)^{b-1}\leq b(b-1) x^{-2} (1-x)^{b-2}\leq 2 x^{-2}\left( 1-(1-x)^b-bx(1-x)^{b-1}\right) , \ x\in [0,1],
\]
which is easy to check, for example, via the  Binomial Theorem.
 
(ii) The stated property of the sequence $\gamma$ was already noted and
used by Schweinsberg, cf. \cite{jS00a} Lemma 3.
For completeness we include a brief argument: From (\ref{Elam_gam})
\begin{eqnarray*}
\gamma_{b+1}-\gamma_b&=&\int_{[0,1]} (x + (1-x)^{b+1} -(1-x)^{b}) x^{-2} d\la(x) \\
&=&\int_{[0,1]} (1-(1-x)^{b}) x^{-1} d\la(x)\geq 0.
\end{eqnarray*}
\hfill \endpf

The following two lemmas and a corollary are auxiliary results, often 
implicitly observed in \cite{jP99} or \cite{jS00a}, and are of interest to anyone
studying fine properties of $\la$-coalescents.
Fix $a\in (0,1)$.
Let $\la^a$ be the restriction of $\la$ to $[0,a]$, namely 
\[
\la^a([0,x]) = \la([0,a]\cap [0,x]),\ x\in[0,1].
\]
Let 
$\lambda_b^a,\gamma_b^a$
be defined in (\ref{Edeflambk})-(\ref{Edefgam}) using 
$\la^a$ as the underlying measure instead of $\la$.

\begin{Lemma}
\label{Lbounds}
(i) For each fixed $a$, such that $\la^a( (0,1) )>0$,
there exists a constant $C_1=C_1(\la,a)\in (0,\infty)$
such that for all large $b,$
\[
\lambda_b^a\leq \lambda_b\leq C_1 \lambda_b^a.
\]
(ii) There exists an $a<1$ and 
$C_2=C_2(\la,a)\in (0,\infty)$ such that for all large $b,$
\[
\gamma_b^a\leq \gamma_b\leq C_2 \gamma_b^a.
\]
(iii) If $\int_{[0,1]}\frac{1}{x} \,d\la(x) =\infty$, in particular if 
(\ref{Edown}) holds, then for each
fixed $a$, the inequalities in  
(ii) hold with a constant $C_3=C_3(\la,a)\in (0,\infty)$.
\end{Lemma}
{\em Remark.}
For any fixed $\la$
let 
\[
\eta_b:=\sum_{k=2}^b {b \choose k}k\,\lambda_{b,k}.
\] 
Then it is easy to see that 
$\gamma_b/\eta_b \to 1$ as $b\to \infty$,
so statements (ii) and (iii) above extend to the corresponding 
$\eta_b$ and $\eta_b^a$.
\hfill $\diamond$

\noindent
{\em Proof.}
For each $a\in (0,1)$, the first inequalities in both (i) and (ii) are trivial 
consequences of $\la^a$ being the restriction of $\la$, the identities in 
(\ref{Elam_gam}), and the fact that $1-(1-x)^b-bx(1-x)^{b-1}$ and
$bx -1 +(1-x)^b$ are both non-negative on $[0,1]$.

The second inequality in (i) is easy as well, since
$1-(1-x)^b-bx(1-x)^{b-1}$ is bounded by $1,$ which implies  
\begin{equation}
\label{Ehelpin}
\lambda_b\leq \lambda_b^a + \frac{1}{a^2} \la([a,1]).
\end{equation}
Then either $\lambda_b\to \infty,$ in which case (\ref{Ehelpin}) 
implies
$\lambda_b^a\to \infty$ as $b\to \infty$, so that 
for all large $b,$  $\frac{1}{a^2} \la([a,1])\leq \lambda_b^a$,
or $\lambda_b$ stays finite, in which case the upper bound is trivial.

The proof of the second inequality in (ii) is similar.
First note that $bx -1 +(1-x)^b\leq bx$ so that
\begin{equation}
\label{Ehelpin1}
\gamma_b\leq \gamma_b^a + \frac{b}{a} \la([a,1]).
\end{equation}
Now it is easy to see by Lemma \ref{Linit}(ii) that $\gamma_{b+1}
-\gamma_b$ is non-decreasing in $b$ so that $\gamma_b \geq
(\gamma_3-\gamma_2) (b-2)$
for each $b$. 
For $a$ chosen sufficiently close to $1,$
$\frac{1}{a} \la([a,1])<(\gamma_3-\gamma_2)/3$ 
(recall $\la$ is a finite measure).
Hence, (\ref{Ehelpin1}) implies 
$\gamma_b^a\geq (\gamma_3-\gamma_2) b/2$ for all $b$ large enough 
and (\ref{Ehelpin1}) then also implies the upper bound in (ii)
since $\frac{b}{a}\la([a,1])<(\gamma_3-\gamma_2)\frac{b}{3}<\gamma^{a}_{b}.$

Part (iii) follows immediately
from the argument for (ii), and the following fact
(already noticed by Pitman \cite{jP99}, Lemma 25),
\begin{equation}
\label{Egamasy}
\int_{[0,1]}\frac{1}{x}\, d\la(x) = \lim_{b\to \infty} \frac{\gamma_b}{b}.
\end{equation}
In particular, (\ref{Edown}) must imply that the left hand side in 
(\ref{Egamasy}) is infinite.
\hfill \endpf

Let the symbol $\asymp$ stand for  
``asymptotically equivalent behavior`` 
in the sense that $a_m \asymp b_m$ (as $m\to \infty$)
if there exist two finite positive constants $c,C$ such that
\[
c \,a_m \leq b_m \leq C \,a_m, \ m\geq 1.
\]
\begin{Lemma}
\label{Lasymptotic}
We have
\begin{itemize}
\item[(i)]
\[
\lambda_b\asymp b^2\la[0,1/b] + \int_{[1/b,1]}\frac{1}{x^2} d\la(x)
- b \int_{[1/b,\log(2(b-1)/(1-e^{-1}))/b)]}\frac{(1-x)^{(b-1)}}{x} d\la(x),
\]
\item[(ii)]
\[
\gamma_b\asymp b^2\la[0,1/b] + b \int_{[1/b,1]}\frac{1}{x} d\la(x).
\]
\end{itemize}
\end{Lemma}
{\em Proof.}
(i) 
To show the first claim, 
use expression (\ref{Elam_gam}) to get for $b\geq 2,$
\begin{eqnarray*}
\lambda_b&=& \int_{[0,1/b]} \frac{1-(1-x)^b-bx(1-x)^{b-1}}{x^2}d\la(x)\\
& &+ \int_{[1/b,1]} \frac{1-(1-x)^b-bx(1-x)^{b-1}}{x^2}d\la(x).
\end{eqnarray*}
Then note that
\[
\frac{\int_{[0,1/b]} \frac{1-(1-x)^b-bx(1-x)^{b-1}}{x^2}d\la(x)}
{b(b-1)\la([0,1/b])/2} \to 1, \mbox{ as } b\to\infty,
\]
and also that
\[
(1-e^{-1})
\int_{[1/b,1]} \frac{1}{x^2} \, d\la(x)
\leq
\int_{[1/b,1]} \frac{1-(1-x)^b}{x^2}\,d\la(x)
\leq
\int_{[1/b,1]} \frac{1}{x^2} \, d\la(x).
\]
A calculus fact, $1-x\leq e^{-x}$, $x\in [0,1]$
implies
that if $x\geq \log(2(b-1)/(1-e^{-1}))/b$, then 
$(1-x)^{b-1}\leq \frac{(1-e^{-1})}{2b}$.
This in turn implies that
\begin{eqnarray*}
& &\int_{[\log(2(b-1)/(1-e^{-1}))/b,1]} \frac{b(1-x)^{b-1}}{x}\, d\la(x)\\
&\leq&
\int_{[\log(2(b-1)/(1-e^{-1}))/b,1]}\frac{(1-e^{-1})}{2x}\, d\la(x)\\
&\leq& 
\frac{(1-e^{-1})}{2}\int_{[\log(2(b-1)/(1-e^{-1}))/b,1]}\frac{1}{x^2}\,
d\la(x)\\
&\leq& \frac{1}{2} \int_{[1/b,1]}\frac{1-(1-x)^b}{x^2}\,
d\la(x),
\end{eqnarray*}
so that $- \int_{[\log(2(b-1)/(1-e^{-1}))/b,1]} \frac{b(1-x)^{b-1}}{x}\,
d\la(x)$
can be ignored in the asymptotics, and the remaining term
\[
\int_{[1/b,\log(2(b-1)/(1-e^{-1}))/b]} \frac{b(1-x)^{b-1}}{x}\,d\la(x),
\]
appears in the asymptotic expression for $\lambda_b$. 
(ii) Since $\gamma_b\asymp \eta_b$, see the above remark,
it suffices to show the second statement for $\eta_b$ instead.
As in (\ref{Elam_gam}),
\[
\eta_b =b\int_{[0,1/b]} \frac{1-(1-x)^{b-1}}{x} \, d\la(x) +
b \int_{[1/b,1]} \frac{1-(1-x)^{b-1}}{x}\, d\la(x),
\]
and since it is easy to see that
\[
\frac{\int_{[0,1/b]} \frac{1-(1-x)^{b-1}}{x} \, d\la(x)}
{b \la([0,1/b])} \to 1, \mbox{ as } b \to \infty,
\]
while
\[
\int_{[1/b,1]} \frac{1-(1-x)^{b-1}}{x}\, d\la(x) \asymp
\int_{[1/b,1]} \frac{1}{x} \, d\la(x),
\]
the claim on the asymptotics of $\gamma_b$ (i.e., $\eta_b$) follows.
\hfill\endpf

\begin{Corollary}
\label{Ceasyasy}
\begin{itemize}
\item[(i)] If $\lambda_b\to \infty$, as $b\to \infty$ then 
$\lim_{b\to \infty} \lambda_{b+1}/\lambda_b = 1$,\\
\item[(ii)]
Since $\gamma_b\to \infty$, as $b\to \infty$ we obtain that 
$\lim_{b\to \infty} \gamma_{b+1}/\gamma_b = 1$. 
\end{itemize}
\end{Corollary}
{\em Proof.}
(i) 
By the Binomial Formula, for $x\in [0,1],$
\[
\frac{b-1}{2} b x^2 (1-x)^{b-1} \leq
1- (1-x)^b - bx(1-x)^{b-1},
\]
so that 
\begin{eqnarray}
\label{Elamcomp}
& &\int_{[0,\log(2(b-1))/(b-1)]} b (1-x)^{b-1} \, d\la(x) \\
\nonumber
&\leq& 
\frac{2}{b-1}
\int_{[0,\log(2(b-1))/(b-1)]}
 \frac{1-(1-x)^b-bx(1-x)^{b-1}}{x^2}d\la(x).
\end{eqnarray}
Since $\int_{[\log(2(b-1))/(b-1),1]} b (1-x)^{b-1} \, d\la(x) \leq 
 \frac{b}{2(b-1)}\la([0,1])$,
the conclusion follows by Lemma \ref{Linit}(i), 
(\ref{Elam_gam}), (\ref{Elamcomp}) and the
fact that $\lambda_b\to \infty$.\\
(ii) 
Perhaps the easiest way to see that $\gamma_b\to \infty$
whenever $\la([0,1])>0$ is by using the identity (\ref{Egamasy}).
The statement then follows immediately
from Lemma \ref{Linit}(i),  Lemma \ref{Lasymptotic}(ii), and the fact that 
\[
\int_{[0,1]} \frac{1-(1-x)^b}{x}  d\la(x) \leq 2 b\la([0,1/b]) +
\int_{[1/b,1]}\frac{1}{x} \,d\la(x). 
\]
\hfill\endpf

\begin{Lemma}
\label{Lpowerlaw}
There exists a finite number $\rho\geq 1$
such that for any $\la$, and all $b,m\geq 2$ such that $b/m\geq 2$
we have
\[
\lambda_{b} \leq m^\rho \lambda_{\lceil b/m\rceil}
\]
\end{Lemma}
{\em Proof.}
In this lemma we consider the identities (\ref{Elam_gam}) for
all real $b\geq 1$.
It suffices to show that 
\begin{equation}
\label{Efirstbound}
\lambda_b \leq c \lambda_{b/2}
\end{equation}
for all $b\geq b_0$ where $b_0$ is some finite integer.
Indeed, if $m\in (2^k, 2^{k+1}]$ for some $k$
then 
\[
\lambda_b\leq c^{k+1} \lambda_{b/2^{k+1}}
\leq c m^{\log_2 c} \lambda_{\lceil b/m\rceil},
\]
and now one can take $\rho > \log_2 c+ \frac{\log c}{\log 2}$
to get the statement of the lemma.
Define the function 
\[
g(\beta,x):= 1-(1-x)^\beta - \beta x(1-x)^{\beta-1}.
\]  
Due to representation (\ref{Elam_gam}) for $\lambda_b$ it then suffices to study
\[
f_b(x) :=\frac{1-(1-x)^b - bx (1-x)^{b-1}}{1-(1-x)^{b/2}-\frac{b}{2} x(1-x)^{\frac{b}{2}}}=\frac{g(b,x)}{g(b/2,x)},
\]
and show
\[
\sup_{x\in [0,1]} f_b(x) \leq c,
\]
uniformly in all $b\geq b_0$.
Note that $f_b(0+) =4$ and that $f_b(1) =1$.
The derivative $f_b'(x)$ can be written as a ratio $fn_b(x)/fd_b(x)$ where
$fd_b(x)\geq 0$ and where
$fn_b(x)$ equals
\[
x(1-x)^{b/2-2} [ b(b-1) (1-x)^{b/2}
- (b/2-1)\frac{b}{2} - 
\frac{b}{2}(\frac{3b}{2}-1)(1-x)^b - b (\frac{b}{2})^2 x(1-x)^{b-1}].
\]
Therefore $fn_b(x)<0$ whenever $b(b-1) (1-x)^{b/2} < (\frac{b}{2}-1)\frac{b}{2}$
and in particular whenever $x>\frac{2}{b}\log 8$ for all $b\geq 4$.
So it suffices to show that 
\[
\sup_{x\in [0,\frac{2}{b}\log{8}]} f_b(x) \leq c.
\]
For this note that $g(b/2,x) \geq {\frac{b}{2} \choose 2} x^2 (1-x)^{\frac{b}{2}-2}$
for any $x\in [0,1]$, and that (by expanding the binomial terms
and noting $x/(1-x) \leq 4 \log{8}$ 
whenever $x <\frac{2}{b}\log{8}$ and $b\geq 10$)
\[
\frac{g(b,x)}{{b\choose 2} x^2 (1-x)^{b-2}}=\sum_{l=2}^b 
\frac{2(b-2)!}{l!(b-l)!}\left(\frac{x}{1-x}\right)^{l-2} \leq 2\cdot 8^4. 
\]
\hfill\endpf

Now we turn to the spatial setting.
Recall that
the vertex set of $\GG$ is $\{g_1,\ldots,g_\mg\}$.
Denote by $\lambda(b_1,b_2,\ldots,b_\mg)$ the total rate of coalescence for the 
configuration with $b_i$ blocks at site $g_i$,
\[
\lambda(b_1,b_2,\ldots,b_\mg)
:=\sum_{i=1}^\mg \sum_{k=2}^{b_i}{b_i \choose k} \lambda_{b_i,k}=\sum_{i=1}^\mg \lambda_{b_i}.
\]
Similarly, let
\[
\gamma(b_1,b_2,\ldots,b_\mg):=\sum_{i=1}^\mg \gamma_{b_i}.
\]
Denote by $\lfloor x \rfloor$ the integer part of the real number $x$ and let
$\lceil x \rceil:=-\lfloor -x \rfloor$.

The following two lemmas will be useful for the proof of the 
characterization result given in Proposition \ref{PzoE}.
\begin{Lemma} 
\label{Lspat_estim}
For all $\mg \geq 1, b_{i} \geq 0, i=1,\dots, \mg$ integers with $\sum_{i=1}^\mg b_{i} >\mg,$
\begin{itemize}
\item[(i)]   $\gamma_{\sum_{i=1}^\mg b_i} \geq \gamma(b_1,b_2,\ldots,b_\mg)
\geq \mg\gamma_{\lfloor \sum_{i=1}^\mg b_i/\mg\rfloor},$
\item[(ii)]  $\mg^{1+\rho}\lambda_{\lceil \sum_{i=1}^\mg b_i/\mg\rceil} \geq \lambda(b_1,b_2,\ldots,b_\mg)
\geq \lambda_{\lceil\sum_{i=1}^\mg b_i/\mg\rceil}.$
\end{itemize}
\end{Lemma}
{\em Proof.}
(i) 
In order to verify the first inequality  we observe that for $x \in [0,1],$
\begin{equation}
\label{EJena}
\mg-1 \geq (\sum_{i=1}^{\mg} (1-x)^{b_{i}}) - (1-x)^{\sum_{i=1}^{\mg} b_{i}}
\end{equation}
since one can simply check that equality holds for $x=0$ and that 
$x \mapsto (\sum_{i=1}^{\mg} (1-x)^{b_{i}}) - (1-x)^{\sum_{i=1}^{\mg} b_{i}}$
is a decreasing function on $[0,1].$ 
Inequality (\ref{EJena}) implies that
\begin{equation*} 
(\sum_{i=1}^{\mg} b_{i})x-1+(1-x)^{\sum_{i=1}^{\mg} b_{i}}
\geq \sum_{i=1}^{\mg} ( b_{i}x-1+(1-x)^{b_{i}} )
\end{equation*}
for all $x \in [0,1].$ The first inequality in (i) now follows from this and 
from (\ref{Elam_gam}), since  
\begin{equation}
\label{Egamlast}
 \gamma(b_1,b_2,\ldots,b_\mg) =\sum_{i=1}^\mg \gamma_{b_{i}}
= 
\int_{[0,1]} x^{-2} \sum_{i=1}^\mg (b_i x-1+(1-x)^{b_i})d\la(x),
\end{equation}
for $b_{i} \geq 0$ (if we set $0^{0}=1$).
The second inequality of (i) 
is immediate if $\sum_{i=1}^{\mg} b_{i}<2\mg.$ Otherwise, we note that
\begin{equation}
\label{EJen}
\sum_{i=1}^\mg(1-x)^{b_i} \geq \mg (1-x)^{ \sum_{i=1}^\mg b_i/\mg}
\end{equation}
by Jensen's Inequality 
since the function $y\mapsto a^y$ is convex for every $a>0$.
Therefore, (\ref{Egamlast}) is bounded below by
\[
\mg \int_{[0,1]} x^{-2} \left ( \beta x - 1+(1-x)^{\beta }\right) d\la(x),
\]
where $\beta = \sum_{i=1}^\mg b_i/\mg$. If $\beta$ is an integer then the last
expression is just $\mg \gamma_{\beta}.$ Now note that 
the function $\beta\mapsto \beta x - 1+(1-x)^{\beta}$
is increasing (for $\beta \geq1$) and this implies the second
inequality in (i).

(ii) 
Use Lemma \ref{Lpowerlaw} to conclude
\[
\lambda(b_1,b_2,\ldots,b_\mg)\leq \mg^{1+\rho} 
\lambda_{\lceil \sum_{i=1}^{\mg} b_i/\mg\rceil}.
\]
The second inequality of (ii) is a simple consequence of the fact that there exists
a $1\leq j \leq \mg$ such that $b_j \geq \lceil \sum_{i=1}^\mg b_i/
\mg\rceil$ 
and Lemma \ref{Linit}(i).
\hfill \endpf

Now consider the coalescent $(\Pi^{\ell}_{n\mg}(t))_{t\geq 0}$ such that its
initial configuration $\Pi^{\ell}_{n\mg}(0)$ has $n$ blocks at each site of $\GG$.
Let 
\begin{equation}
\label{Tndef}
T_n:=\inf\{t>0: \# \Pi_{n\mg}(t) \leq 2\mg\}. 
\end{equation}
\begin{Lemma}
\label{LETn}
If condition (\ref{Edown}) holds then
$\sup_n E[T_n] \leq \sum_{b=2}^{\infty} \frac{3 \mg^{\rho+1}}{\gamma_{b}} <\infty.$
\end{Lemma}
{\em Proof.}
The argument is an adaptation of the argument by Schweinsberg 
\cite{jS00a}, Lemma 6, to 
our situation.
In fact we will even use similar notation.
For $n\in {\bf N}$ define $R_0:=0$ and stopping times 
(with respect to the filtration generated by $\Pi^{\ell}_{n\mg}$) given by
\begin{eqnarray*}
R_i&:=&1_{\{\# \Pi_{n\mg}(R_{i-1})> 2\mg \}}\inf\{t>R_{i-1}: \#\Pi_{n\mg}(t)<\# \Pi_{n\mg}(R_{i-1})\}\\
& &+ 1_{\{\# \Pi_{n\mg}(R_{i-1})\leq 2\mg \}} R_{i-1} 
,\ i\geq 1. 
\end{eqnarray*}
In words, $R_i$ is the time of the $i$th coalescence as long as the number of blocks
before this coalescence exceeds $2\mg,$ otherwise $R_i$ is set equal to the 
previous coalescence time.
Since there are no more than $2\mg$ blocks left after $(n-2)\mg$ coalescence events,
note that 
\[
T_n =R_{(n-2)\mg}.
\]
Of course, it is also possible that 
$T_n = R_i$ for $i< (n-2)\mg$, but the above identity 
holds almost surely as $R_{(n-2)\mg}=R_i$ in this case.
Let
\[
L_i=R_i-R_{i-1}, \ J_i= \# \Pi_{n\mg}(R_{i-1}) - \# \Pi_{n\mg}(R_i),
\]
and note that there exists some finite random number $\xi_i$ such that
$R_{i-1}=T_0^i<T_1^i<T_2^i<\ldots<T_{\xi_i}^i<R_i$,
where $T_1^i,T_2^i,\ldots$ are the successive times of 
migration jumps of various blocks
from site to site in between the $i-1$th and $i$th coalescence time.
Let $B_i(t)$ be the number of blocks located at site $g_i \in {\cal G}$ at 
time $t.$
Since the total number of blocks does not change at the jump times
$T^{i}_{j}$ for $j=1,\dots,\xi_i$ we have due to Lemma 
\ref{Lspat_estim} (ii) that $\lambda(B_1(T^{i}_{j}),\ldots,B_\mg(T^{i}_{j}))
\geq \lambda_{\lceil\sum_{e=1}^\mg B_e(T^{i}_{j})/\mg\rceil}
=\lambda_{\lceil\sum_{e=1}^\mg B_e(R_{i-1})/\mg\rceil}.$ 
This
implies (by coupling of exponentials in a straightforward way)  that 
\begin{equation}
\label{ELi}
E[L_i|\Pi_{n\mg}(R_{i-1})] \leq 
\frac{1}{\lambda_{\lceil\sum_{e=1}^\mg B_e(R_{i-1})/\mg\rceil}}.
\end{equation}
Also note that for all $i$ with $\Pi_{n\mg}(R_{i-1})>2\mg$,
\begin{eqnarray}
E[J_i|\Pi_{n\mg}(R_{i-1})]\bck&=&\bck
E\left[\left.\sum_{j=0}^\infty\frac{\gamma(B_1(T_j^i),
\ldots,B_\mg(T_j^i))}{\lambda(B_1(T_j^i),\ldots,B_\mg(T_j^i))} 
1_{\{\xi_i=j\}}\right|\Pi_{n\mg}(R_{i-1})\right]
\nonumber\\
&\geq &\frac{\mg\gamma_{\lfloor\sum_{e=1}^\mg B_e(R_{i-1})/\mg\rfloor}}{
 \mg^{1+\rho} \lambda_{\lceil\sum_{e=1}^\mg B_e(R_{i-1})/\mg\rceil}}
P(1_{\{\xi_i<\infty\}}|\Pi_{n\mg}(R_{i-1}))\nonumber\\
&= &\label{EJi}
\frac{1}{\mg^\rho}\frac{\gamma_{\lfloor\sum_{e=1}^\mg B_e(R_{i-1})/\mg\rfloor}}{
\lambda_{\lceil\sum_{e=1}^\mg B_e(R_{i-1})/\mg\rceil}},
\end{eqnarray}
where the first equality is a direct consequence of definitions 
(\ref{Edeflam}) and (\ref{Edefgam}), and the fact that $J_i$ is the
decrease in the number of blocks at the $i$th coalescence time $R_i$.
The middle inequality is due to Lemma \ref{Lspat_estim} (i) and (ii). 
From (\ref{ELi}) and (\ref{EJi}) and the fact that $L_{i}=0$ if $J_{i}=0$
we get the important relation
\begin{equation}
\label{Eimpo}
E[L_i|\Pi_{n\mg}(R_{i-1})]\leq \frac{\mg^\rho}{\gamma_{\lfloor\sum_{e=1}^\mg B_e(R_{i-1})/\mg\rfloor}} 
E[J_i|\Pi_{n\mg}(R_{i-1})]
\end{equation}
for $i \geq 1.$
Now
\begin{eqnarray*}
\bck\bck E[T_n] \bck &=& \bck E[\sum_{i=1}^{\mg(n-2)} L_i]=\sum_{i=1}^{\mg(n-2)} 
E\left[E\left[L_i|\Pi_{n\mg}(R_{i-1})\right]\right]\\
\bck &\leq&\bck  \sum_{i=1}^{\mg(n-2)} 
E\left[ \frac{\mg^\rho}{\gamma_{\lfloor\sum_{e=1}^\mg B_e(R_{i-1})/\mg \rfloor}} 
E\left[J_i|\Pi_{n\mg}(R_{i-1})\right]\right]\\
\bck &=&\bck  E\left[\sum_{i=1}^{\mg(n-2)}  
\frac{\mg^{\rho}}{\gamma_{\lfloor\sum_{e=1}^\mg B_e(R_{i-1})/\mg \rfloor}} J_i\right]=
E\left[\sum_{i=1}^{\mg(n-2)} 
\sum_{j=0}^{J_i-1} \frac{\mg^{\rho}}{\gamma_{\lfloor\sum_{e=1}^\mg B_e(R_{i-1})/\mg \rfloor}} \right]\\
\bck &\leq&
 \mg^{\rho} E\left[ \sum_{b=2}^{n} \frac{\mg}{\gamma_b} +
\frac{2\mg}{\gamma_2} \right]
\leq \sum_{b=2}^{n} \frac{3\mg^{\rho+1}}{\gamma_b},
\end{eqnarray*}
where we have used 
Lemma \ref{Ldetermini} below.
\hfill \endpf
\begin{Lemma}
\label{Ldetermini}
For a fixed $\mg$, let $m,n$ be positive integers such that
$m\in [n\mg,(n+1)\mg)$.  
For any $k\geq 1$ and $j_1,\ldots,j_k\geq 1$ such that
$\sum_{i=1}^{k-1} j_i<m - 2 \mg$ and 
$\sum_{i=1}^k j_i\in[m-2\mg,m-1]$ one has
\begin{equation}
\label{Edeterm}
\sum_{i=1}^k \frac{j_i}{\gamma_{\lfloor (m-\sum_{\ell=1}^{i-1} j_\ell)/\mg
\rfloor}}
\leq 
\frac{m-n\mg}{\gamma_n} + \sum_{b=2}^{n-1}\frac{\mg}{\gamma_b} + 
\frac{2\mg}{\gamma_2}.
\end{equation}
\end{Lemma}
{\em Proof.}
Statement (\ref{Edeterm}) can be proved for each fixed $\mg$
by induction in $n$.
The base cases $n=2$ with $m>2\mg$ and $\sum_{i=1}^{k} j_i=m-1$ explain
 the extra summands $2\mg/\gamma_2$.
Here one also uses the fact that $(\gamma_{b})_{b=2}^{\infty}$ is an increasing sequence
(cf. Lemma \ref{Linit} (ii)). 
\hfill \endpf\\

Let us now recall the construction in Theorem \ref{Tconstr},
and define
\[
T_n^{(2)} := T_n \mbox{ from definition (\ref{Tndef}), }
\]
and
\[
T_\infty = T_{\infty}^{(2)}
= \sup_n T_n^{(2)} =\inf \{ t > 0 : \#\Pi (t) \leq 2\mg\},
\]
and furthermore define
\begin{equation}
\label{Tndefk}
T_n^{(k)}:=\inf\{t>0: \# \Pi_{n\mg}(t) \leq k\mg\}, \
T_\infty^{(k)}:=\sup_n T_n^{(k)},\ k \geq 3. 
\end{equation}
Note that by monotone convergence $T_n^{(k)} \nearrow T_\infty^{(k)}$ we have
\[
E[T_{\infty}^{(k)}] = \lim_{n\to \infty} E[ T_n^{(k)}], \ k\geq 2.
\]
\begin{Corollary}
\label{Cgeneralk}
If condition (\ref{Edown}) holds then for each $k\geq 2,$
$\sup_n E[T_n^{(k)}] \leq \sum_{b=k}^{\infty} \frac{
\mg^{\rho+1}}{\gamma_{b}} + \frac{ k\mg^{\rho+1}}{\gamma_k} <\infty$,
and in particular
\[
\lim_{k\to \infty} \sup_n E[T_n^{(k)}]=0.
\]
\end{Corollary}
{\em Proof.}
The upper bound on $E[T_n^{(k)}]$  
can be shown as in the proof of Lemma \ref{LETn}. 
The second claim above now follows by relation (\ref{Egamasy}) and the
observation following it.
\hfill \endpf\\

We can now establish the following analogues to
Proposition 23 of Pitman \cite{jP99} and 
Proposition 5 of Schweinsberg \cite{jS00a} in 
the spatial setting.
\begin{Proposition}
\label{PzoE}
Assume that $\la$ has no atom at 1. Then the $\la$-coalescent
either comes down from infinity or it stays infinite. Furthermore,
it stays infinite if and only if $E[T_{\infty}]=\infty$.  
\end{Proposition} 
{\em Proof.}
Define
$T:=\inf\{ t\geq 0:  \#\Pi(t) < \infty\}$.
The first statement could be shown
following Pitman \cite{jP99} Proposition 23
by observing that 
$P[ 0< T< \infty]>0$ leads to a contradiction. 
We choose a different approach, based on Corollary \ref{Cgeneralk} and
coupling with non-spatial coalescents.

Suppose that (\ref{Edown}) holds.
Then $E[T_{\infty}]< \infty$, by Lemma \ref{LETn},
 implying $T_{\infty}<\infty$ almost surely.
Also note that the $\la$-coalescent comes down from infinity due to Corollary 
\ref{Cgeneralk}, since for any $t>0$, and any $k\geq 2,$
\[
P[T>t] \leq P[T_\infty^{(k)}>t] \leq \frac{E[T_\infty^{(k)}]}{t}.
\]
This verifies that $P[T \in \{0,\infty\}] =P[T=0] =1,$ again by Corollary \ref{Cgeneralk}.

If (\ref{Edown}) does not hold, we will show next by a coupling 
argument that, provided $\#\Pi(0)=\infty$, we have
$P[T \in \{0,\infty\}]= P[T = \infty]=1$.
This implies of course that
$P[T_\infty = \infty] =1$ and $E[T_\infty]=\infty$.
So assume that $\#\Pi(0)=\infty,$ i.e.~that there exists at least
one site $g$ in $\GG$ such that $\Pi^\ell(0)$ contains infinitely many blocks with 
label $g.$
Then the spatial coalescent $\Pi^{\ell}$ is stochastically bounded below by a
coalescing system $\tilde{\Pi}^{\ell}$, in which 
any block that attempts to migrate is assigned to a ``cemetery site''
$\partial$ instead. 
More precisely, the evolution of the process $\tilde{\Pi}^{\ell}$ at each site is independent 
from the evolution at any other site, and its transition mechanism is
specified by: 
\begin{itemize}
\item[(i)] blocks coalesce according to a $\la$-coalescent,
\item[(ii')] each block vanishes (moves to $\partial$) at rate 1.
\end{itemize}
By adapting the construction of $\Pi^{\ell}$ in Theorem
\ref{Tconstr},
one can easily construct a coupling $(\Pi^{\ell}(t),\tilde{\Pi}^{\ell}(t))_{t\geq 0}$
on the same probability space, 
so that at each time $t$, and for
each site $g$ of $\GG$,
the number of blocks in  $\Pi^{\ell}(t)$  located at $g$ is larger than (or equal to) the
number of blocks in $\tilde{\Pi}^{\ell}(t)$ located at $g$.
We will show that in any given 
time interval $[0,t]$, at each site of $\GG$ that initially 
contained infinitely many blocks,  
there are infinitely
many blocks remaining in $\tilde{\Pi}^{\ell}$ (even though there are infinitely many blocks
that do vanish to $\partial$ by time $t$). 
Therefore, $\infty = \#\tilde{\Pi}(t)\leq \#\Pi(t)$ so that $\Pi^{\ell}$
stays infinite.

To show that $P[\#\tilde{\Pi}(t)=\infty]=1$ for each $t>0$, it will 
be convenient to construct a
coupling of $\tilde{\Pi}^{\ell}(t)$ with a new random object $\Pi^1(t)$. 
Since there is no interaction between the sites of $\GG$ in $\tilde{\Pi}$,
it suffices to consider the nonspatial case where $|\GG|=1$.
Introduce an auxiliary family 
$(X_j)_{j\geq 1}$ of independent exponential random variables
with parameter $1$. 
Take a (non-spatial) $\la$-coalescent $(\Pi^1(s))_{s \in[0,t]}$ 
such that $\Pi^1(0)= \tilde{\Pi}(0)$,
and in addition augment the state space for $\Pi^1$ to accommodate a
mark for each block. 
Initially all blocks start with an empty mark.
At any $s\leq t$, any block $A\in \Pi^1(s)$ is marked by
$\partial$ if $\{X_{\min A} \leq s\}$.
In this way, if an already marked block $A$ coalesces with a family $A_1,A_2,\ldots$
of blocks, such that $\min A\leq \min_j (\min A_j)$, the new block
$A\cup \cup_j A_j$
automatically inherits the mark $\partial$.
Note as well that if a marked block $A$ coalesces with at least one unmarked block
containing  a smaller element than $\min A$, the new block will be
unmarked.

The number $\#^u\Pi^1(t)$ of all {\em unmarked} blocks in $\Pi^1(t)$ is
stochastically smaller than the number $\# \tilde{\Pi}(t)$.
To see this, note the difference between $\tilde{\Pi}(t)$ and $\Pi^1(t)$:
a marked block in $\Pi^1(s)$ is not removed from the population
immediately (unlike in $\tilde{\Pi}$) so it may coalesce (and ``gather'') additional
blocks with higher indexed elements during $[s,t]$ resulting in
a smaller number of unmarked partition elements in $\Pi^1(0)$ than in $\tilde{\Pi}(t)$.

Another random object $\Pi^2(t)$, equal in distribution 
to $\Pi^1(t)$, can be constructed as follows:
run a (non-spatial) $\la$-coalescent
$(\Pi^2(s))_{s\in [0,t]}$, and attach to each block $A\in \Pi^2(t)$ a mark 
$\partial$ with probability $e^{-t}$.
Let $\#^u\Pi^2(t)$ be the number of all unmarked blocks in $\Pi^2(t)$. 
Since (\ref{Edown}) does not hold,
due to the corresponding result in \cite{jS00a},
$P[\#\Pi^2(t)=\infty]=1.$
Since $\Pi^1(t)$ and $\Pi^2(t)$ have the same distribution,
then $\#^u \Pi^1(t)$ and $\#^u \Pi^2$
have the same distribution and 
by the above construction we conclude immediately that
\[
1=P[\#^u\Pi^2(t)=\infty]=P[\#^u\Pi^1(t)=\infty].
\]
Recalling that $\#^u\Pi^1(t)$ is stochastically bounded above by $\tilde{\Pi}(t)$ for all $t\geq 0$
completes the proof.
\hfill \endpf\\

{\em Remark.}
It is intuitively clear that in the case in which the $\la$-coalescent $\Pi^{\ell}$ stays infinite,
there are infinitely many  blocks in $\Pi^{\ell}$ at all positive times at 
all sites, a proof of this fact is left to an interested reader.
\hfill $\diamond$

\section{Uniform asymptotics}
\label{UNI}
Note that the upper bound in Lemma
\ref{LETn} and Corollary \ref{Cgeneralk} 
neither depends on the structure of $\GG$
nor on the underlying migration mechanism.
After a careful look at the proofs the reader will see that 
in fact the same estimates would hold with an arbitrary migration
mechanism, even if it is not independent 
from the coalescent mechanism.

In this section we will use the fact that each block 
changes its label (i.e.~migrates) at rate $1$, independently
from the coalescent mechanism.
Recall the setting of Lemma \ref{LETn} and Corollary 
\ref{Cgeneralk}.
\begin{Theorem}
\label{Tuniform}
If  (\ref{Edown}) holds,
then there exists a constant $c$ uniform in $\la$, $\mg$, 
the structure of $\GG$, and the transition kernel of the
migration mechanism, such that 
\[
\sup_{n} E [T_{n}] \leq  \sum_{b=2}^{\infty} \frac{1}{\gamma_{b}}
+ \frac{2}{\gamma_{2}},
\]
and moreover
\[
\sup_{n} E [T_{n}^{(k)}] \leq 
\left(\sum_{b=k}^{\infty} \frac{1}{\gamma_{b}}+\frac{k}{\gamma_k}\right).
\]
\end{Theorem} 
{\em Proof.}
We use the same notation as in the proof of Lemma \ref{LETn}, 
but this time the calculations are finer.
First, fix an $i \geq 1$ (note the
subscripts $i$ are omitted in a number of places below for notational convenience).
Recall the jump times 
$T_{j}^{i}$ and configurations $\Pi_{n\mg}(T^i_j)$ with 
\[
a:=\sum_{e=1}^\mg B_e(R_{i-1})
= \sum_{e=1}^\mg B_e(T_{j}^{i})
\mbox{ for all }j\leq \xi_i.
\] 
Also set
$\lambda_j^i=\lambda_{j}=\lambda(B_1(T_j^i),\ldots,B_\mg(T_j^i))$ for all 
$j\in \bN_{0}$.
Recall that $\xi_i:=\max\{k: T_k^i<R_i\}$ is the number of migration events
in between $(i-1)$st and $i$th coalescence time.
Note that the quantities $\lambda_j$ are relevant for our
process only if $j \leq \xi_{i}$.

Using the first line of (\ref{EJi}) and Lemma \ref{Lspat_estim} (i) 
as well as further conditioning on $(\lambda_l)_{l \in \bN_{0}}$
we obtain
\begin{eqnarray}
E[J_i|\Pi_{n\mg}(R_{i-1})]\bck&\geq &\bck
\mg \gamma_{\lfloor\frac{a}{\mg}\rfloor}
E\left[\left.\sum_{j=0}^\infty \lambda_{j}^{-1}
1_{\{\xi_i=j\}}\right|\Pi_{n\mg}(R_{i-1})\right]
\nonumber\\
&= & \mg \gamma_{\lfloor \frac{a}{\mg} \rfloor}
E\left[\left.\sum_{j=0}^\infty \lambda_{j}^{-1}
P\left[\xi_i=j | (\lambda_l)_{l \in \bN_{0}},\Pi_{n\mg}(R_{i-1})\right]   
\right|\Pi_{n\mg}(R_{i-1})\right]
\nonumber\\
\label{Jieq}
&= & \mg \gamma_{\lfloor \frac{a}{\mg} \rfloor}
E\left[\left.\sum_{j=0}^\infty \lambda_{j}^{-1}
\left( \prod_{l=0}^{j-1} \frac{a}{\lambda_{l} +a} \right) 
\frac{\lambda_{j}}{\lambda_{j} +a} \right|\Pi_{n\mg}(R_{i-1})\right].
\end{eqnarray}
For the next computation define an auxiliary 
i.i.d.~sequence $(X_{j})_{j \geq 0}$  of exponential random variables 
with parameter $a$, as well as
a sequence $(Y_{j})_{j \geq 0}$ of independent random variables 
where each $Y_j$ has an exponential ( $\lambda_{j}$) distribution.
Note that
$W_{j} := X_{j} \wedge Y_{j}$ are exponential random variables with rate $a+ \lambda_{j}$ 
that are independent from $Z_{j}= 1_{\{X_{j}>Y_{j}\}}.$ 

Observe that conditioned on $((\lambda_{j})_{j \in \bN_{0}}, \Pi_{n\mg}(R_{i-1}))$ 
the $X_{j}$ correspond to
the waiting time until the next migration and the $Y_{j}$ to the waiting
time until coalescence as long as $\sum_{l=1}^{j-1} Z_{l}=0.$ So the event 
$\{Z_{0}=\dots=Z_{j-1}=0 \} = \{ \xi_{i}\geq j\}$ is independent of $W_{j}$.
This implies that
\begin{eqnarray}
\nonumber
E\left[L_i | \Pi_{n\mg}(R_{i-1})\right] 
&=& E\left[ \left.E\left[ \left.\sum_{j=0}^{\xi_i} W_{j} 
 \right| (\lambda_{l})_{l \in \bN_{0}},\Pi_{n\mg}(R_{i-1}) \right] \right|\Pi_{n\mg}(R_{i-1})\right]\\
\nonumber
&=& E\left[ \left.E\left[ \left.\sum_{j=0}^{\infty} 
W_{j}  1_{\{\xi_{i}\geq j\}}
 \right| (\lambda_{l})_{l \in \bN_{0}},\Pi_{n\mg}(R_{i-1}) \right]
\right|\Pi_{n\mg}(R_{i-1})\right]\\
\nonumber
&=& E \left[ \sum_{j=0}^{\infty}  
E \left[ \left. W_j \right| 
(\lambda_{l})_{l \in \bN_{0}},\Pi_{n\mg}(R_{i-1}) \right] \right.
\\
\nonumber
& & \quad \quad \cdot  E\left[ \left. 1_{\{\xi_{i}\geq j\}}
\right| (\lambda_{l})_{l \in \bN_{0}},\Pi_{n\mg}(R_{i-1}) \right]
 \bigg|\ \Pi_{n\mg}(R_{i-1})\Bigg]\\
\label{Lieq}
&=& E\left[ \left. \sum_{j=0}^{\infty}  
\frac{1}{\lambda_j+a} 
\left( \prod_{l=0}^{j-1} \frac{a}{\lambda_{l}+a}\right)
 \right|\Pi_{n\mg}(R_{i-1})\right].\\
\end{eqnarray}

Comparing now the terms in (\ref{Jieq}) and (\ref{Lieq}) 
 we find that 
\begin{eqnarray}
\label{Eimpo2}
E\left[L_i | \Pi_{n\mg}(R_{i-1})\right] \leq
\frac{1}{\mg \gamma_{\lfloor\frac{a}{\mg}\rfloor}} E[J_i|\Pi_{n\mg}(R_{i-1})] ,  
\end{eqnarray}
where we gained a factor of $\mg^{\rho+1}$ in the denominator with respect to
the analogous relation
(\ref{Eimpo}) in the proof of Lemma \ref{LETn}. The rest of the proof
proceeds now as the proof of Lemma \ref{LETn} and Corollary \ref{Cgeneralk}
and hence we obtain
\begin{equation*}
E[T_{n}] \leq  \sum_{b=2}^{n} \frac{1}{\gamma_b} + \frac{2}{\gamma_2},
\end{equation*}
and 
\[
E[T_{n}^{(k)}] \leq \left(\sum_{b=k}^{n} \frac{1}{\gamma_b} + 
\frac{k}{\gamma_k}\right),
\]
\hfill \endpf\\
{\bf Definition.}
We will say that the $\la$-coalescent {\em comes down from infinity uniformly}
if 
\[
\lim_{k\to \infty} \sup_n E T_n^{(k)} = 0.
\]
\hfill $\diamond$

In particular, by Proposition \ref{PzoE} and Theorem 
\ref{Tuniform} any coalescent with independent Markovian migration mechanism
that comes down from infinity also
comes down from infinity uniformly.\\
{\bf Example.}
Let $\alpha\in (0,2)$.
The Beta($2-\alpha,\alpha$)-coalescent, where
$\la\ed{\rm Beta}(2-\alpha,\alpha)$
has density $x^{1-\alpha}(1-x)^{\alpha-1}/\Gamma(2-\alpha)\Gamma(\alpha)$
is of special 
interest in \cite{BBCEMSW05}.
As already noted in \cite{jS00a}, for $\alpha \in (0,1]$ this
(non-spatial) coalescent stays infinite, and for $\alpha\in (1,2)$ it
comes down from infinity.
By the previous theorem the spatial  
Beta($2-\alpha,\alpha$)-coalescent 
comes down   from infinity uniformly. 
An interesting consequence follows
 by the results of the next section.
\hfill $\diamond$

\section{Asymptotics on large tori}
\label{ASY}
In this section we further restrict the setting in the following way:
\begin{itemize}
\item
the graph $\GG$ is a $d$-dimensional torus $T^N=[-N,N]^d \cap \bZ^d$ for some $N\in
\bN$, where $d\geq 3$ is fixed,
\item
the migration corresponds to a random walk on the torus, meaning that
the kernel $p(x,y), x,y\in \GG$ is given
as $p(x,y)\equiv \sum_{\{z: (z-y)\!\! \mod N=0\}}\tilde{p}(z-x)$,
where $\tilde{p}$ is purely $d$-dimensional
distribution such that $\sum_x |x|^{d+2} \tilde{p}(x) <\infty$,
\item
the $\la$-coalescent comes down  from infinity (uniformly), i.e.,
condition (\ref{Edown}) holds.
\end{itemize}

We are concerned here with convergence of the $\la$-coalescent partition
structure  on $T^{N},$ if time is rescaled by the volume $(2N+1)^{d}$ 
of $T^N$,
to that of a time-changed {\em non-spatial} Kingman coalescent as $N
\rightarrow \infty.$
The main results are presented in Theorem \ref{TnpartCon} and Theorem
\ref{TuniCon}:
Theorem  \ref{TnpartCon} states convergence of the partition structure in
a functional sense
for arbitrary finite initial configurations. Theorem \ref{TuniCon} states
convergence of the number of
partition elements in a functional sense if the initial number of
partition elements is infinite.

We write $\PP^{N,\ell}$ if we want to emphasize that  the partitions are labeled by $T^N$.
Let 
\[
\Pi^{N,\ell}_{\pi} \
\mbox{ and } \
\Pi^{N,\ell}
\]
denote the $\la$-coalescent started from a
partition $\pi \in \PP^{N,\ell},$ and  the
$\la$-coalescent started from
any partition that contains 
infinitely many equivalence classes 
labeled by (located at) each site of
$T^N$, respectively. In order to determine the large space-time asymptotics for
$\Pi^{N,\ell}$, at time scales on the order of the volume $(2N+1)^d$ of
$T^N,$
we imitate a ``bootstrapping'' argument from \cite{GLW}.

{\em Remark.}
Observe that in \cite{GLW}, only the singular $\la=\delta_0$ case
was studied in this context.
However, the structure of the argument
concerning large space-time asymptotics
carries over due to the cascading property for
general (spatial) $\la$-coalescents, in particular due to the fact that
any two partition elements $\pi_1,\pi_2\in \Pi^{N,\ell}(0)$
coalesce
at rate
\begin{equation}
\label{Ega2}
\lambda_{2,2} = \la([0,1])
\end{equation}\
while they are at the same site, and that they do not coalesce otherwise.
\hfill $\diamond$\\

We will need the following notation: for a marked partition $\pi \in
\PP_n^\ell$ (or $\pi \in \PP^\ell$),
and two real numbers $a<b\in \bR$, write
\[
\pi \in [[a,b]],
\]
if $\forall i ,j$ with $ i\neq j$, such that
$(A_i,\zeta_i),(A_j,\zeta_j)\in \pi$
we have $|\zeta_i - \zeta_j| \in [a,b]$.
In words, $\pi \in [[a,b]]$ if and only if
all the mutual distances for pairs of different partition
elements of $\pi$ are contained in $[a,b]$.

The following theorem states that, viewed on the right timescale
$t(2N+1)^{d},$ and after some initial collapse of a finite starting
configuration, the partitions of the $\la$-coalescent on the tori $T^N$
with $N$ large behave like those of a (non-spatial) time-changed 
Kingman coalescent. 
To make this statement more precise, we introduce the following notation.

Let $G=\sum_{k =0}^{\infty}\tilde{p}_{k}(0)$ where $\tilde{p}_{k}$ denotes the $k$-step
transition probability of a $\tilde{p}$ random walk.
Note that this random walk is transient on $\bZ^d$, so that $G<\infty$.
Let $\Pi_\pi^{\bZ^d, \ell}$ be the $\la$-coalescent on $\GG=\bZ^{d}$ with
migration given by the random walk kernel $\tilde{p}$,
started from partition $\pi$ with $\#\pi<\infty$.
The transience of $\tilde{p}$ also implies existence of
non-trivial limit partitions
\[
\Pi_{\pi}^{\bZ^d}(\infty)=\lim_{t \rightarrow
\infty} \Pi_{\pi}^{\bZ^d}(t),
\]
in the sense that if $\#\pi \geq 2$ then
$\#\Pi_{\pi}^{\bZ^d}(\infty)\geq 2$
with positive probability.

We define $K_{\pi}$ as the {\em non-spatial} Kingman coalescent
started in the partition $\pi \in \PP$ or $\pi\in \PP_n.$ This means that
$K_{\pi}$ is the $\la_{K}$-coalescent for $\la_{K}= \delta_{0}$ and
$|\GG|=1$ with initial configuration $K_{\pi}(0)=\pi.$

Denote by $D(\bR_+, E)$
 the c\`adl\`ag paths on $\bR_+$ with  values in some metric space
$E$, and equip the space
$D(\bR_+, E)$ with the usual Skorokhod
topology. Also let $"\Rightarrow"$ indicate convergence in distribution.
Set
\begin{equation}
\label{Ekappa}
\kappa= \frac{2}{G+ 2/\lambda_{2,2}}.
\end{equation}
Recall that $\Pi^{N,\ell}$ starts from a configuration containing infinitely
many blocks, namely  the partition $\Pi(0) = \{\{j\}: j \in \bN\}$.
The theorem below concerns the behavior of only finitely many
blocks.
Recall that $\Pi^{N,\ell}(0)|_n$ is the restriction of the labeled
partition
$\Pi^{N,\ell}(0)$ to $[n]$.
In the theorem below
we use the abbreviation
$\Pi^{N,\ell}_n:= \Pi^{N,\ell}_{\Pi^{N,\ell}(0)}|_n$.
Again, $\Pi^N_n$ is the process of
partitions corresponding to $\Pi^{N,\ell}_n$.

\begin{Theorem}
\label{TnpartCon}
Assume that 
 for each fixed $n\geq 1$ and 
all large $N$ we have $\Pi^{N,\ell}(0)|_n =\Pi^{N+1,\ell}(0)|_n$. 
Then for each $n$,  we obtain as
$N \rightarrow \infty,$
the following convergence 
of the (unlabeled) partition processes:
\[
(\Pi_{n}^{N}(t (2N+1)^{d}))_{t \geq 0} \Rightarrow
(K_{\Pi_{n}^{\bZ^d}(\infty)}
(\kappa t) )_{t \geq 0},
\]
where convergence is with respect to
the Skorokhod topology
on $D(\bR_+, \PP_n),$ and both $\Pi_{n}^{N,\ell}$ and
$\Pi_{n}^{\bZ^d,\ell}$ are started from the same initial configuration
$\Pi^{N,\ell}(0)|_n \in \PP_{n}^{\ell}.$
\end{Theorem}

{\em Remark.}
The statement is a generalization of  Proposition 7.2 in \cite{GLW}, which
deals with the case of spatial Kingman coalescents, rather than
$\la$-coalescents and only states convergence of the marginals.
Nevertheless, the first part of the argument is analogous, and  we will
change it only slightly in preparation for Proposition \ref{PmarginCon} and 
Theorem \ref{TuniCon}.
\hfill $\diamond$\\
As the first step we will state a result for the case in which the initial
configuration is
sparse on the torus, so that no coalescence involving more than two
particles may be
seen in the limit. The general case, stated in Theorem \ref{TnpartCon},
will then follow easily.
\begin{Proposition}
\label{PnpartCon}
Let $a_{N}\rightarrow \infty$ be such that $a_{N}/N \rightarrow 0.$
Fix $n \in \bN$, and let
$\pi^{N,\ell} \in \PP^{N,\ell}$ be such that
$\#\pi^{N,\ell}=n\geq 2$,
$\pi^{N} \in [[a_{N},\sqrt{d}N ]]$,
and such that its corresponding (unlabeled) partition
$\pi^N$ equals a constant partition $\pi_0\in \PP$ for all $N$.
Then as $N \rightarrow \infty$, we have the following convergence in
distribution of the (unlabeled) partition processes:
\[
(\Pi_{\pi^{N,\ell}}^{N}(t (2N+1)^{d}))_{t \geq 0} \Rightarrow
(K_{\pi_0}
(\kappa t) )_{t \geq 0},
\]
where the convergence is in the space $D(\bR_+, \PP).$
\end{Proposition}
{\em Proof.}
To simplify the notation we refer to the $i$th block of $\pi_0$ as
$\{i\}$, for $i=1,\ldots n$.
In order to show the convergence on the space $D(\bR_+, \PP)$
we will prove that the joint distribution of
inter-coalescence  times converges, when appropriately rescaled, to the joint
distribution of inter-coalescence times of $K(\kappa \,\cdot)$,
and that, at each coalescence time, any pair of remaining blocks is equally
likely to
coalesce next, see also \cite{GLW2} for a similar argument.

We set $\tau_{0}^{N}=0.$ Since there are at most $n-1$ coalescence times
in general,
we then define recursively stopping times for $k=1,\dots,n-1,$
\[
\tau_{k}^{N}:=\inf\{t \geq \tau_{k-1}^{N}
: \#\Pi_{\pi^{N,\ell}}^N(t)\neq \#\Pi_{\pi^{N,\ell}}^N(\tau_{k-1}^{N})\},
\]
as long as $\#\Pi_{\pi^{N,\ell}}^N(\tau_{k-1}^{N})>1.$ 
Also define inter-coalescence times $\sigma_k^N:=
\tau_{k}^N-\tau_{k-1}^N$, $k\leq n-1$.
Let us first observe that 
for $n=2$
\begin{equation}
\label{EfromGLW}
P[\sigma_{1}^{N}/(2N+1)^{d}<t]  =P[\tau_{1}^{N}/(2N+1)^{d}<t]  \to e^{-\kappa t}
\end{equation}
uniformly in $t\in[0,T]$ for any $T<\infty$,
by Lemma 7.3 in \cite{GLW}. 
Indeed, as remarked at the beginning of this section, the 
spatial $\la$-coalescent restricted to two-particles is identical in
law to the spatial $\lambda_{2,2} \delta_0(\cdot)$-coalescent 
from \cite{GLW}.

Let $U_{k}$ be independent
exponential random variables with parameters
$\kappa \binom{n-(k-1)}{2}$ for $k<n-1.$
We wish to show the convergence in distribution of the random vector
\begin{equation}
\label{taujointconv}
(\sigma_{1}^{N} /(2N+1)^{d}, \dots, \sigma_{n-1}^{N}/(2N+1)^{d}) \Rightarrow
(U_{1}, \dots, U_{n-1})
\end{equation}
as $N \rightarrow \infty.$
The statement is clear by (\ref{EfromGLW}) if $n=2$.
 In order to show (\ref{taujointconv}) for $n \geq 2,$
the first step is to see that, 
we may exclude the possibility of coalescence of more than
two particles at any given time with probability tending to $1$ as $N
\rightarrow \infty$.

Let $\tau^{N}(i,j)$ be the time of the coalescence which merges the block
$A^{(i)}$
containing $i$ and the block $A^{(j)}$ 
containing $j$, and for each $i$ denote by $\zeta^{(i)}$
the label associated
with the block $A^{(i)}$.
Then, 
we have
for any $0<T<\infty,$ and any distinct $i,j,k \in [n],$ 
\begin{equation}
\label{nothree}
\int_{0}^{T(2N+1)^{d}}
P\left[ \tau_{1}^{N}=\tau^{N}(i,j) \in du, |\zeta^{(i)}-\zeta^{(k)}|\leq
a_{N}  \right]
\to 0,
\end{equation}
uniformly over all
partitions $\pi^{N,\ell}\in[[a_N,\sqrt{d}N]]$,
as $N \to \infty$.
The statement (\ref{nothree})
is analogous to (3.7) in Cox \cite{Cox89},
and follows with exactly the same
calculation. 
Likewise, a statement analogous to (3.8) in \cite{Cox89} holds, saying that 
uniformly over all
$\pi^{N,\ell}\in[[a_N,\sqrt{d}N]]$
\begin{equation}
\label{notclose}
\int_{0}^{T(2N+1)^{d}}
P\left[ \tau_{1}^{N}=\tau^{N}(i,j) \in du, |\zeta^{(k)}-\zeta^{(l)}|\leq
a_{N}  \right]
\to 0,
\end{equation}
as $N \to \infty$ for $i,j,k,l \in [n]$ distinct. 

Now fix $T<\infty$, $\epsilon>0$, and let $n>2$.
Relation (\ref{nothree}) implies that for $N$ large
enough,
\begin{equation}
\label{1less}
P\left[ \# \Pi_{\pi^{N,\ell}}^N(\tau_{1}^{N})\neq n-1,\tau_{1}^{N}< T(2N+1)^{d}\right]<\epsilon,
\end{equation}
and together with (\ref{notclose}) it implies
that 
\[
P\left[ \Pi_{\pi^{N,\ell}}^{N,\ell}(\tau_{1}^{N})\not\in [[a_N, \sqrt{d}N]],\tau_{1}^{N}<
T(2N+1)^{d}\right]<\epsilon.
\]
A simple induction (using the strong Markov property and uniformity of
(\ref{nothree}) and (\ref{notclose}) in $t\in [0,T]$)
yields the following statement: 
for each $k<n-1$, and any fixed $\eps>0$, 
if $N$ is large enough then
\[
P\left[ \# \Pi_{\pi^{N,\ell}}^N(\tau_{k}^{N})\neq n-k,\tau_{k}^{N}< T(2N+1)^{d}\right]<\epsilon,
\]
and 
\[
P\left[ \Pi_{\pi^{N,\ell}}^{N,\ell}(\tau_{k}^{N})\not\in [[a_N, \sqrt{d}N]],\tau_{k}^{N}<
T(2N+1)^{d}\right]<\epsilon,
\]
for $k\leq n-2$.
From this we get that, for any fixed $\eps>0$, 
if $N$ is large enough, 
\begin{equation}
\label{Eimpo1}
P\left[\# \Pi_{\pi^{N,\ell}}^N(\tau_{k}^{N})= n-k 
\mbox{ for each $k$ with } \tau_{k}^{N}<T(2N+1)^d\right]> 1-\eps, 
\end{equation}
and
\begin{equation}
\label{Eimpo2a}
P\left[\Pi_{\pi^{N,\ell}}^N(\tau_{k}^{N})\in [[a_N, \sqrt{d}N]] 
\mbox{ for each $k$ with } \tau_{k}^{N}<
T(2N+1)^{d}\right]>1-\epsilon.
\end{equation}
Moreover,
on the event 
\[
\{\tau_{k}^{N}<
T(2N+1)^{d}\}\cap \{ \# \Pi_{\pi^{N,\ell}}^N(\tau_{k}^{N}) = n-k\} 
\cap \{ \Pi_{\pi^{N,\ell}}^N(\tau_{k}^{N})\in [[a_N, \sqrt{d}N]]\}
\]
we have as in (3.1) of \cite{Cox89} that
\begin{equation}
\label{Eimportant}
|P[\sigma_{k+1}^N/ (2N+1)^d > u|\FF_{\tau_{k}^{N}}] 
- e^{-\kappa {n-k \choose 2} u }|<\eps_N,
\end{equation}
where $\eps_N$ depends on $N$ only, and where $\eps_N \to 0$,
as $N\to \infty$.

In order to arrive at (\ref{taujointconv}),
we show that  $\sigma_{k}^{N}$ is 
asymptotically independent of $\sigma_{k-1}^{N}, \ldots, \sigma_{1}^{N}$ for
all $k=2, \dots, n-1.$ So consider
for any fixed $0 \leq t_{1},\dots,t_{k}$, where $\sum_{i=1}^k t_i <T,$ the event
\[
A_{k}^{N}:=\left\{\frac{\sigma_{k}^{N}}{(2N+1)^d}< t_{k}, 
\frac{\sigma_{k-1}^{N}}{(2N+1)^d}< t_{k-1}, \dots, 
\frac{\sigma_{1}^{N}}{(2N+1)^d}<
t_{1} \right\}.
\]
In particular, on this event we have that $\tau_{i}^{N}< T(2N+1)^d$ is satisfied for $i=1,\dots, k.$ 
We obtain
\begin{eqnarray*}
& &P\left[  A_{k}^{N} \right]\\
&=&E\left[ P\left[\left.\frac{\sigma_{k}^{N}}{(2N+1)^d}<t_{k}\right|{\cal F}_{\tau_{k-1}^{N}}
\right]
1_{ A_{k-1}^{N} }\right]\\
&=&
E\left[ \left(
P\left[\left.\frac{\sigma_{k}^{N}}{(2N+1)^d}<t_{k}
\right|{\cal F}_{\tau_{k-1}^{N}}\right]
- \left(1-e^{-\binom{n-(k-1)}{2}\kappa t_{k}}\right)
\right)
1_{ A_{k-1}^{N}}\right]\\
& &+(1-e^{-\binom{n-(k-1)}{2}\kappa t_{k}})
P\left[ A_{k-1}^{N}\right].
\end{eqnarray*}
Now use (\ref{Eimpo1}), (\ref{Eimpo2a}), and (\ref{Eimportant}) 
to get 
\[
\lim_{N\to \infty}
P\left[ A_{k}^{N}  \right] = 
(1-e^{-\binom{n-(k-1)}{2}\kappa t_{k}}) 
\lim_{N\to \infty}
P\left[  A_{k-1}^{N} \right].
\]
By iterating the argument we obtain asymptotic independence.
This in turn implies that $(\#\Pi_{\pi^{N,\ell}}^N(t
(2N+1)^{d})_{t \geq 0}
\Rightarrow (\# K_{\pi_{0}}(\kappa t))_{t \geq 0}$ in the Skorokhod topology,
since by (\ref{Eimpo1}), as $N \to \infty,$
\[
P\left[\#\Pi_{\pi^{N,\ell}}^N(t)=n-\sum_{k=1}^{n-1} 1_{ \{\tau_{k}^{N}<t\}} \text{ for all } 
t < (2N+1)^d T \right] \to 1,
\]
so that  the convergence of the jump times $\tau_{k}^{N}$ in $n-\sum_{k=1}^{n-1} 1_{ \{\tau_{k}^{N}<t\}}$ implies convergence in the Skorokhod topology, see for example Proposition
6.5 in Chapter 3 of \cite{EK}.

Finally, (2.8) in \cite{Cox89} states that for the $\tilde{p}$ random
walk on $T^{N},$
 \[
 \lim_{N \to \infty} \sup_{t\geq (\log N) N^{2}} \sup_{x \in T^{N}}
(2N+1)^{d}|\tilde{p}(x,0)-(2N+1)^{-d}|=0.
 \]
This implies that 
the positions of partition elements in 
$\Pi_{\pi^{N,\ell}}^{N,\ell}(\tau^{N}_{k}+(\log N) N^{2})$
(note that $(\tau^{N}_{k}+(\log N) N^{2})/(2N+1)^{d} \approx
\tau^{N}_{k}/(2N+1)^{d}$)
are approximately uniformly and independently distributed on the
torus.
Due to (\ref{taujointconv}), with probability tending to $1$ 
as $N \rightarrow \infty$, we also have 
\[
\#\Pi_{\pi^{N,\ell}}^N(\tau^{N}_{k}+(\log N) N^{2} ) =
\#\Pi_{\pi^{N,\ell}}^N (\tau^{N}_{k}).
\] 
Therefore, at time $\tau_{k+1}^{N},$ each pair
of partition elements of 
$\#\Pi_{\pi^{N,\ell}}^N (\tau^{N}_{k})$
is approximately 
equally likely to coalesce, as is the case in the
Kingman coalescent. This completes the proof of convergence
on the space $D(\bR_{+}, \PP).$
\hfill \endpf

\bigskip
{\em Proof of Theorem \ref{TnpartCon}.}
Fix $n\in\bN$.
We will first show that, as $N\to \infty$,
 $\Pi_{n}^{N}(N^{3/2})=
\Pi_{n}^{\bZ^d}(\infty)$ (note this is only a statement about
the partition structure, not the locations), and that
$ \Pi_{n}^{N,\ell}(N^{3/2}) \in [[N^{3/4}/\log N,\sqrt{d}N]]$,
with probability arbitrarily close to $1$.
The statement of the theorem 
will then follow by Proposition \ref{PnpartCon} if we 
continue running 
the process from time $N^{3/2}$ 
onwards, and use the strong Markov property, noting that $N^{3/2}=o((2N+1)^d)$.

First define the stopping time
\[
\tau^{N}:=\inf \{ t >0: \max \{\zeta : (A, \zeta) \in \Pi_{n}^{N,\ell}(t)
\} \geq N\}.
\]
Before time $\tau^{N}$ none of the blocks have reached the boundary 
of $[-N,N]^d$, so 
we may couple $\Pi_{n}^{N,\ell}$ and $\Pi_{n}^{\bZ^d, \ell}$ 
in a natural way such that
$\Pi_{n}^{N,\ell}(t)= \Pi_{n}^{\bZ^d,\ell}(t)$ for $t \leq \tau^{N}.$

Note  that by
the functional CLT, any random walk $X$ on $\bZ^d$ with random walk kernel
$\tilde{p}$ started at $X(0)\leq \frac{N}{2}$ satisfies
\[
\lim_{N \rightarrow \infty} P\left[\sup_{0\leq t \leq N^{3/2}
}X(t)<N \right]=1.
\]
Since for $N$ large enough,  $\max \{\zeta : (A, \zeta) \in \Pi^{N,\ell}(0)|_{n} \} \leq
\frac{N}{2}$ and since the coalescent has at most $n$ blocks independently
performing  random walks, we immediately obtain 
\[
\lim_{N \rightarrow \infty} P\left[\tau^N> N^{3/2}\right]=1.
\]
In particular, we have
\begin{equation}
\label{Eintermed}
\lim_{N\to \infty} 
P\left[ \Pi_{n}^{N,\ell}(N^{3/2})=
\Pi_{n}^{\bZ^d,\ell}(N^{3/2})\right]=1.
\end{equation}

To see  that the  blocks remaining at time $N^{3/2}$ are at a
mutual distance of $N^{3/4}/ \log N$ 
with high probability, more precisely
that
\begin{equation}
\label{unispaced}
\lim_{N \rightarrow \infty}P[ \Pi_n^{N,\ell}(N^{3/2}) \in
[[N^{3/4}/\log{N},\sqrt{d}N]] \, ]=1,
\end{equation}
if suffices to observe that again by the functional CLT,
\[
\lim_{N \to \infty}P\left[|X^1(N^{3/2})-X^2(N^{3/2})|<N^{ \frac{3}{4} }/ \log N\right]=0,
\]
where $X^1$ and $X^2$ are two independent $\tilde{p}$-random walks on
$\bZ^{d}$ started at $X^1(0)=X^2(0)=0.$ Due to  (\ref{unispaced}), 
and the fact that
$P[X^1(t)= X^2(t)  \mbox{ for some } t\geq 0 | X^1_0-X^2_0=x]\to 0$ as 
$|x|\to \infty$ we have,
\begin{equation}
\label{toPiinfty}
\lim_{N \rightarrow \infty}P\left[ \Pi_{n}^{\bZ^d,\ell}(N^{3/2})=
\Pi_{n}^{\bZ^d,\ell}(\infty)\right]=1.
\end{equation}
Now (\ref{Eintermed})  and (\ref{toPiinfty}) imply
\begin{equation}
\label{toPiinfty1}
\lim_{N \rightarrow \infty}P\left[ \Pi_{n}^{N,\ell}(N^{3/2})=
\Pi_{n}^{\bZ^d,\ell}(\infty)\right]=1.
\end{equation}
\hfill \endpf

We will also show a uniform convergence to the Kingman coalescent, 
on the same time scale, in the sense of the number of blocks,
cf. Theorem \ref{TuniCon} below.
One starts with a bound on the mean number
of partition elements left in the coalescent $\Pi^N$
at a fixed time, say $1$.
The following useful monotonicity property
carries over from the spatial Kingman coalescent setting to 
the spatial $\la$-coalescent setting:

Suppose that the partition elements of $\Pi^{N,\ell}(0)$ are initially divided into
classes $\Pi^{N,1,\ell}(0),$ $\Pi^{N,2,\ell}(0),\ldots$
(in any prescribed deterministic way)
and let $(\cup_j \Pi^{N,j,\ell}(t))_{t\geq 0}$ denote the united $\la$-coalescent
where only elements of the same class are allowed to coalesce.
\begin{Lemma}
\label{Lcoupl}
For each $t>0,$
\[
E[\#\Pi^N(t)] \leq \sum_j E[\#\Pi^{N,j}(t)].
\]
\end{Lemma}
{\em Proof.}
We can couple $\Pi^{N,\ell}$ and $\cup_j \Pi^{N,j,\ell}$,
using the same Poisson point process (from the construction of
$\Pi^{N,\ell}$)
for all the $\la$-coalescents corresponding to different classes.
It then follows that  $\Pi^N(t)$ is a coarser partition
than $\cup_j \Pi^{N,j}(t)$ for each $t$, almost surely.
This gives the inequality, $\#\Pi^N(t)\leq \sum_j \Pi^{N,j}(t)$,
and in particular the bound in expectation holds.
\hfill \endpf\\
The following lemma is taken from \cite{GLW} and is
similar to Theorem~1 in \cite{BraGri80} and the proposition in
Section 4 of \cite{Cox89}.
\begin{Lemma}
\label{Lnumber}
There is a finite constant $c_d$ such that
uniformly in $N\in\bN$,
and in the sequences $(\Pi^N(0))_{N\in\bN}$ satisfying
$\#\Pi^N(0)\geq (2N+1)^d$,
\begin{equation*}
  {\bf E}\left[\#\Pi^N(t)\right]
 \le
  c_d\max\left\{1,\frac{\#\Pi^N(0)}{t}\right\}.
\end{equation*}
\end{Lemma}
{\em Proof.}
All we need to do is translate the notation and explain
the small differences in the argument.

Our $\lambda_{2,2}$ is $\gamma$ in \cite{GLW}.
The migration walk $\tilde{p}$ is from the same class as in \cite{GLW}.
There are only two
statements in the argument of
\cite{GLW}, Lemmas 7.4 and 7.5 that depend on the structure of
the underlying coalescent.
One is relation (7.50) at the beginning of the argument of Lemma 7.4.
Take $A_0\in \Pi^N(0)$
and note that, similar to (7.44) in \cite{GLW},
\[
\#\Pi^N(t)\leq \#\Pi^N(0) - \sum_{\Pi^N(0)\ni A \neq A_0}
1_{\{A_0 \sim_{\Pi^N(t)} A\}},
\]
so that
\[
E \left[\#\Pi^N(t)\right]\leq E\left[\#\Pi^N(0)\right] - \sum_{\Pi^N(0)\ni
A \neq A_0}
P[A_0 \sim_{\Pi^N(t)} A],
\]
leading to (7.46) of Lemma 7.4 in \cite{GLW}, and therefore to relation (7.50)
since the remaining calculations concern the behavior of two partition
elements
(not the joint behavior of several partition elements).

The other statement concerns (7.58) in the proof of Lemma 7.5: here, the
torus is cut up into boxes and (7.58) states that the expected number of
blocks is bounded by the expected number of blocks in a coalescent in
which only blocks that start in the same initial box may coalesce. This
holds in our setting due to Lemma \ref{Lcoupl}.

Given (7.50) and (7.58), the remaining arguments are the same as those in 
the proof of Lemmas 7.4 and 7.5 of \cite{GLW}.
\hfill \endpf

The next lemma says that the number of the partition elements at
time $\eps (2N+1)^d$ is tight in $N$.
\begin{Lemma}
\label{Ldominating}
Fix $\,0<\epsilon,
\epsilon'<1.$ Then there exists a 
constant $M^0=M^0(\eps,\eps')$ such that, for all 
$M \geq M^{0}$,
\begin{equation*}
\limsup_{N \rightarrow \infty} P[ \#\Pi^{N}(\epsilon (2N+1)^{d}) >M]\leq \epsilon'.
\end{equation*}
\end{Lemma}
{\em Proof.}
Assume $1<\frac{\epsilon (2N+1)^{d}}{2}.$ 
Due to Theorem \ref{Tuniform}, for $k \in \bN,$
\begin{eqnarray*}
\sup_N P[\#\Pi^{N}(1) > k (2N+1)^{d}] &=&\sup_N
P[T^{N,(k)}_{\infty}>1] \\
&=&
\sup_N \sup_n E[ T^{N,(k)}_{n}]
\leq  \left( \sum_{b=k}^{\infty} \frac{1}{\gamma_b} +
\frac{k}{\gamma_{k}}\right),
\end{eqnarray*}
Due to (\ref{Edown}) (more precisely observation (\ref{Egamasy})), the right hand side
converges to zero as
$k \rightarrow \infty.$ Therefore, we may choose 
$M_0\geq 1$ large enough so that 
$c(\sum_{b=M_0}^{\infty} \frac{1}{\gamma_b} +M_0/\gamma_{M_0}) <\eps'/2$ 
and also that 
$M_0> \frac{4c_d}{\epsilon \epsilon'}$.
Then for all $M\geq M_0,$
\begin{equation}
\label{Efirest}
P[\#\Pi^{N}(1) > M (2N+1)^{d}] \leq \frac{\epsilon'}{2}.
\end{equation}
Now take $M\geq M_0$ and
define the event $A^{N}_{M}:=\{ \#\Pi^{N}(1) \leq 
M (2N+1)^{d}\}.$ 
We then have by Lemma \ref{Lnumber} that
\[
E[ \#\Pi^{N}(\epsilon (2N+1)^{d}) | A^{N}_{M} ] \leq c_d \max
\{1,\frac{M (2N+1)^{d}}{\epsilon(2N+1)^{d}-1} \}\leq c_d
\max\{ 1,  \frac{2 M}{\epsilon}\}.
\]
{\em Remark.} Note that on $A_M^N$ we may have $\# \Pi^N(1)\geq (2N+1)^d$ 
and we can apply Lemma \ref{Lnumber} directly, otherwise 
couple the coalescent $(\Pi^N(t), t\geq 1)$ with another coalescent
$\tilde{\Pi}^N(t),t\geq 1)$ such that $\tilde{\Pi}^N$ almost surely
dominates $\Pi^N(t)$ at all times, at all sites, and such that 
$\#\Pi^N(0) = (2N+1)^d$, and apply Lemma \ref{Lnumber} to $\tilde{\Pi}^N$.
\hfill $\diamond$

It follows that
\[
P[\#\Pi^{N}(\epsilon (2N+1)^{d}) >M^2 |A^{N}_{M} ]  \leq
\frac{1}{M} c_d \max\{ 1,  \frac{2}{\epsilon}\}.
\]
By conditioning on whether $A^N_{M}$ or
its complement occurs, using (\ref{Efirest})
\[
P[\#\Pi^{N}(\epsilon (2N+1)^{d}) >M^2 ] \leq \frac{1}{M} c_d \max\left\{ 1, 
\frac{2}{\epsilon}\right\} +
1 \cdot \frac{\epsilon'}{2}.
\]
Since $\frac{2 c_d}{M \epsilon}< \frac{\epsilon'}{2}$ we arrive at
\[
\sup_{N} P[\#\Pi^{N}(\epsilon (2N+1)^{d}) >M^2 ] \leq \epsilon',
\]
for $M\geq M_0$, which gives 
the statement of the lemma with $M^0=(M_0)^2$.
\hfill \endpf

As a consequence, we obtain the following asymptotics for the number of
partitions in $\Pi^N$, a spatial $\la$-coalescent  
started from a partition having infinitely many equivalence classes
labeled by (located at) each site of $T^N$.
\begin{Proposition}
\label{PmarginCon}
Let $(K(t))_{t\geq 0}$ be the (non-spatial) Kingman coalescent started
from the partition $K(0)=\{\{i\},i\in \bN\}$, and let $\kappa$
be defined in (\ref{Ekappa}).
Then, for each fixed $t>0$, we have 
\[
\#\Pi^{N}(t(2N+1)^d)\Rightarrow
\#K(\kappa t),
\]
as $N\to \infty$, where the above convergence is in distribution.
\end{Proposition}
{\em Proof.}
We start with a lower bound on 
the asymptotic distribution of $\#\Pi^{N}(t(2N+1)^d)$.
Let $a_N=N^{3/2}$ so that 
$a_N\to \infty$ and also $\sqrt{a_N}/N \to 0$. 
For any fixed $M$ one can find $N_0$ large enough 
so that for all $N\geq N_0$, $\Pi^{N,\ell}(0)$ contains at least $M$
blocks (say $A_{i_1},\ldots, A_{i_M}$),
having mutual distances larger than $\sqrt{a_N}$.
Let $\tilde{\Pi}^{N,\ell}$ be the $\Pi^{N,\ell}$ coalescent restricted to 
$\{A_{i_1},\ldots, A_{i_M}\}$.
Then clearly 
\begin{equation}
\label{Eupbdhelp}
\#\Pi^N(t(2N+1)^d) \geq \#\tilde{\Pi}^N(t(2N+1)^d).
\end{equation}
As a consequence of
(the proof of) Theorem \ref{TnpartCon}, for any $t>0$ , as $N\to \infty,$
\[
P[\#\tilde{\Pi}^N(t(2N+1)^d) =k ] \to P[\#K_M(\kappa t) = k], \ k=1,\ldots,M,
\]
where $K_M(\cdot)$ is the Kingman coalescent 
started from partition $\{\{1\},\ldots,\{M\}\}$.
By (\ref{Eupbdhelp}), for $k =1,\ldots,M$,
\begin{eqnarray*}
\liminf_{N\to \infty} P[ \#\Pi^N(t(2N+1)^d) \geq k]&\geq&
 \liminf_{N\to \infty} P[ \#\tilde{\Pi}^N(t(2N+1)^d) \geq k]\\
&=& P[K_M(\kappa t) \geq k].
\end{eqnarray*}
Taking $M\to \infty$ on both sides and using the well-known 
coming down (or entrance law)
property for $K(\cdot)$, we get for each $k\geq 1$, and each $t>0$,
\begin{equation}
\label{Elobd}
\liminf_{N\to \infty} P[ \#\Pi^N(t(2N+1)^d) \geq k]\geq
P[K(\kappa t) \geq k].
\end{equation}
Before continuing, note an interesting consequence:
If $t_N = o ((2N+1)^d)$ then
\begin{equation}
\label{Etclzero}
\lim_{N \to \infty}P[\#\Pi^N(t_N) \geq k]\geq \lim_{t\searrow 0} \liminf_{N\to \infty} P[ \#\Pi^N(t(2N+1)^d) \geq
k]=1, \ k\geq 1,
\end{equation}
or equivalently, $\#\Pi^N(t_N) \to \infty$ in probability as $N\to \infty$.
To get the upper bound corresponding to (\ref{Elobd}), we use 
Lemma \ref{Ldominating}. Namely, fix $\eps,\eps'\in (0,1\wedge t/2)$, 
and find the corresponding $M^0=M^0(\eps/2,\eps')$.
Running the configuration $\Pi^{N,\ell}(\eps(2N+1)^d/2)$ for an additional 
$N^{3/2}<<\eps(2N+1)^d/2$ units of time will result in 
$\Pi^{N,\ell}(\eps(2N+1)^d/2+N^{3/2})$.
On the event $\{\#\Pi^{N}(\eps(2N+1)^d/2)\leq M^0\}$, that has
probability greater than $1-\eps'$, we have 
$\{\#\Pi^{N}(\eps(2N+1)^d/2+N^{3/2})\leq M^0\}$, and moreover
due to (\ref{unispaced}), for $N$ sufficiently large, 
all the (fewer than $M^0$) partition elements of 
$\Pi^{N,\ell}(\eps(2N+1)^d/2+N^{3/2})$ are at mutual distances larger
than $N^{3/4}/\log{N}$ 
with probability close to $1$.
More precisely, if we let $C_{N,\eps,M^0}$  be the event that
\begin{eqnarray*}
&&\Pi^{N,\ell}(\eps(2N+1)^d/2+N^{3/2}) \in [[N^{3/4}/\log{N},\sqrt{d}N]] \text{ and }\\
&&\#\Pi^N(\eps(2N+1)^d/2+N^{3/2})\leq M^0
\end{eqnarray*}
then, for all sufficiently large $N$,
\begin{equation}
\label{Enew}
P[C_{N,\eps,M^0}] \geq 1-2\eps'.
\end{equation}
Again by the proof of Theorem \ref{TnpartCon}, 
on $C_{N,\eps,M^0}$ 
we have for $k=1,\ldots,M^0$, almost surely
\[
|P[\#\Pi^N(t(2N+1)^d)\geq k|\FF_{\eps(2N+1)^d/2+N^{3/2}}]
- P[\#K_{\#\Pi^N(\eps(2N+1)^d/2+N^{3/2})}(\kappa t)\geq k]|\leq \eps_N,
\]
further implying, 
\begin{equation}
\label{Esortofbd}
|P[\#\Pi^N(t(2N+1)^d)\geq k, C_{N,\eps,M^0}] -
E[P[\#K_{\#\Pi^N(\eps(2N+1)^d/2+N^{3/2})}(\kappa t)\geq k]1_{C_{N,\eps,M^0}}]|\leq \eps_N,
\end{equation}
where $\eps_N \to 0$ as $N\to \infty$.
Now use 
\begin{eqnarray*}
P[\#K_{\#\Pi^N(\eps(2N+1)^d/2+N^{3/2})}(\kappa t)\geq k] 1_{C_{N,\eps,M^0}}
&\leq&
P[\#K_{M^0}(\kappa t)\geq k] 1_{C_{N,\eps,M^0}}\\
&\leq&
P[\#K(\kappa t)\geq k]1_{C_{N,\eps,M^0}}
\end{eqnarray*}
together with the fact that $C_{N,\eps,M^0}^c$ happens with
probability smaller than $2\eps'$ to obtain from (\ref{Esortofbd})
that 
\[
\limsup_{N\to \infty} 
P[ \#\Pi^N(t(2N+1)^d) \geq k]
\leq 4\eps' + P[\#K(\kappa t)\geq k].
\]
The last statement is true for any $\eps'>0$, and this combined with
(\ref{Elobd}) gives
\begin{equation}
\label{Elim}
\lim_{N\to \infty} P( \#\Pi^N(t(2N+1)^d) \geq k)=
P(K(\kappa t) \geq k), \ k \geq 1.
\end{equation}
\hfill \endpf\\

An even stronger form of convergence is true.
It holds in any setting where Proposition \ref{PmarginCon}
and Theorem \ref{TnpartCon} hold, in particular in the setting of 
\cite{GLW}, although there it does not appear explicitly.
Its analogue is important for the diffusive clustering analysis 
in the two-dimensional setting of \cite{GLW2}.
\begin{Theorem}
\label{TuniCon}
Let $(K(t))_{t\geq 0}$ be 
as in Proposition \ref{PmarginCon}.
Then for each fixed $a>0$, we have 
\[
(\#\Pi^{N}(t(2N+1)^d))_{t\geq a}\Rightarrow
(\#K(\kappa t))_{t\geq a},
\]
as $N\to \infty$, where the convergence is with respect to
the Skorokhod topology on c\`adl\`ag processes.
\end{Theorem}
{\em Proof.}
As a consequence of (\ref{Enew}) we have for any fixed $a>0$, 
\begin{equation}
\label{apart}
\lim_{N\to \infty} P[ \Pi^{N,\ell}(a(2N+1)^d) \in [[N^{3/4}/\log{N},\sqrt{d}N]]]= 1.
\end{equation}
Together with the convergence of marginals in Proposition 
\ref{PmarginCon}, and  Theorem \ref{TnpartCon}, this 
yields the current statement.
\hfill\endpf\\

\noindent
{\bf Acknowledgement.} 
We thank Anita Winter and Robin Pemantle for
very useful discussions, and Ed Perkins
for careful reading of a preliminary draft. A.S.~would like to thank 
the Mathematics Department of the
University of British Columbia for its hospitality.

\end{document}